\documentclass[journal]{new-aiaa}
\usepackage[utf8]{inputenc}

\usepackage{graphicx}
\usepackage{amsmath}
\usepackage[version=4]{mhchem}
\usepackage{siunitx}
\usepackage{longtable,tabularx}
\usepackage{indentfirst}
\usepackage{hyperref}
\usepackage{mathtools}
\usepackage{bm}
\usepackage{xspace}
\usepackage{wrapfig}
\usepackage{multicol}
\usepackage{xurl}
\usepackage[flushleft]{threeparttable}
\usepackage{booktabs}
\usepackage{gensymb}
\usepackage{mdframed}
\usepackage{subcaption}
\usepackage{xcolor, pifont}
\usepackage{lscape}
\usepackage{tcolorbox}
\usepackage[ruled, vlined, linesnumbered]{algorithm2e}

\let\svthefootnote\thefootnote
\newcommand\freefootnote[1]{
  \let\thefootnote\relax
  \footnotetext{#1}
  \let\thefootnote\svthefootnote
}

\newcommand{\EOSSP}{\hyperlink{EOSSP}{\textsf{EOSSP}}\xspace}
\newcommand{\REOSSP}{\hyperlink{REOSSP}{\textsf{REOSSP}}\xspace}
\newcommand{\REOSSPExact}{\hyperlink{REOSSP}{\textsf{REOSSP}}\textsf{-Exact}\xspace}
\newcommand{\RHP}{\hyperlink{RHP}{\textsf{REOSSP-RHP}}\xspace}
\newcommand{\qedsymbol}{\hfill $\square$}

\definecolor{myblue}{rgb}{0, 0.23, 0.64}
\definecolor{WVUblue}{rgb}{0, 0.16, 0.33}
 
\hypersetup{
    colorlinks=true,
    linkcolor=myblue,
    filecolor=magenta,
    citecolor=myblue,
    urlcolor=myblue,
}

\title{Reconfigurable Earth Observation Satellite Scheduling Problem}

\author{Brycen D. Pearl\footnote{Ph.D. Student, Department of Mechanical, Materials and Aerospace Engineering, Student Member AIAA.}, Joseph M. Miller\footnote{Undergraduate Student, Department of Mechanical, Materials and Aerospace Engineering.}, and Hang Woon Lee\footnote{Assistant Professor, Department of Mechanical, Materials and Aerospace Engineering; hangwoon.lee@mail.wvu.edu. Member AIAA (Corresponding Author).}}
\affil{West Virginia University, Morgantown, WV, 26506}

\begin{document}

\newpage

\freefootnote{This paper is a substantially revised version of the paper AIAA 2025-0589, presented at the 2025 AIAA SciTech Forum, Orlando, FL, January 6-10, 2025. It offers new results, an additional solution methodology, and a better description of the materials.}

\maketitle

\begin{abstract}
    Earth observation satellites (EOSs) play a pivotal role in capturing and analyzing planetary phenomena, ranging from natural disasters to societal development. The EOS scheduling problem (EOSSP), which optimizes the schedule of EOSs, is often solved with respect to nadir-directional EOS systems, thus restricting the observation time of targets and, consequently, the effectiveness of each EOS. This paper leverages state-of-the-art constellation reconfigurability to develop the reconfigurable EOS scheduling problem (REOSSP), wherein EOSs are assumed to be maneuverable, forming a more optimal constellation configuration at multiple opportunities during a schedule. This paper develops a novel mixed-integer linear programming formulation for the REOSSP to optimally solve the scheduling problem for given parameters. Additionally, since the REOSSP can be computationally expensive for large-scale problems, a rolling horizon procedure (RHP) solution method is developed. The performance of the REOSSP is benchmarked against the EOSSP, which serves as a baseline, through a set of random instances where problem characteristics are varied and a case study in which Hurricane Sandy is used to demonstrate realistic performance. These experiments demonstrate the value of constellation reconfigurability in its application to the EOSSP, yielding solutions that improve performance, while the RHP enhances computational runtime for large-scale REOSSP instances.
\end{abstract}

\section*{Nomenclature}
{
\renewcommand\arraystretch{1.0}
\noindent\begin{longtable*}{@{}l @{\quad=\quad} l@{}}
$B_{\text{various}}$ & battery-related values \\
$C$                  & arbitrary weight of data downlink in objective functions \\
$c$                  & cost of orbital maneuver \\
$c_{\max}$           & budget for total orbital maneuver costs \\
$D_{\text{various}}$ & data-related values \\
$G$                  & number of ground stations \\
$\mathcal{G}$        & set of ground stations \\
$H$                  & binary visibility condition of the sun \\
$J$                  & number of orbital slot options \\
$\mathcal{J}$        & set of orbital slot options \\
$K$                  & number of satellites \\
$\mathcal{K}$        & set of satellites \\
$L$                  & number of lookahead stages \\
$\mathcal{L}$        & set of lookahead stages \\
$P$                  & number of targets for observation \\
$\mathcal{P}$        & set of targets \\
$S$                  & number of reconfiguration stages \\
$\mathcal{S}$        & set of reconfiguration stages \\
$T$                  & number of time steps \\
$T_r$                & finite schedule duration \\
$\mathcal{T}$        & set of time steps \\
$V$                  & binary visibility condition of targets \\
$W$                  & binary visibility condition of ground stations \\
$Z$                  & figure of merit \\
$z$                  & objective function value \\
$\Delta t$           & time step size between discrete time steps \\
$b$                  & indicator variable of the current battery storage level \\
$d$                  & indicator variable of the current data storage level \\
$g$                  & ground station index \\
$h$                  & decision variable to control solar charging \\
$i, j$               & orbital slot option indices \\
$k$                  & satellite index \\
$\ell$               & lookahead stage index \\
$p$                  & target index \\
$q$                  & decision variable to control data downlink to ground stations \\
$s$                  & reconfiguration stage index \\
$t$                  & time step index \\
$x$                  & decision variable to control constellation reconfiguration \\
$y$                  & decision variable to control target observation \\
\end{longtable*}
}

\section{Introduction}

\lettrine{S}{atellite} systems are crucial for gathering information on various planetary phenomena on Earth through Earth observation (EO), leveraging remote sensing to provide observations in the form of radio frequencies, radar, lidar, and optical imaging, among other measurements. EO data may be used to benefit studies of direct societal activity, such as in use observing agricultural drought progression \cite{Crocetti2020Drought}, providing additional security and situational awareness through maritime surveillance \cite{Soldi2021Maritime}, and providing insight into current and potential civil infrastructure development \cite{Prakash2020Urban}. Additionally, EO data may be used to provide necessary information regarding natural processes, such as when utilized to monitor and report on natural disasters \cite{Chien2019Taskable} and in use modeling the potential spread of diseases \cite{Parselia2019Disease}. 

Optimally collecting and providing EO data requires considering the \textit{Earth observation satellite scheduling problem} (EOSSP), an optimization problem where each satellite task is treated as a decision variable to maximize observation rewards. Typically, formulations of the EOSSP include the scheduling of observation tasks to gather data, data downlink tasks to transmit data to ground stations, and physical or operational constraints \cite{Cho2018Optimization,Herrmann2023Reinforcement,Eddy2021}. An important constraint of many EOSSP formulations is the restriction of observation tasks only to times when targets are visible to satellites \cite{CHATTERJEE2022}, defined as the \textit{visible time window} (VTW) \cite{Chen2018MILP}. As a result, the amount of data collected through observation is directly influenced by the VTW, further underlining the VTW as one of the most impactful aspects of any EOSSP. 

The VTW, and consequently the maximum amount of collectible data, are maintained as fixed parameters due to the underlying assumptions present in EOSSP formulations. Two such underlying assumptions are the restriction of the pointing direction of satellites, employing only nadir-directional observations that present a fixed VTW, and the restriction to a given orbit that is incapable of being manipulated \cite{Herrmann2022}. Such limitations imply that the VTW is a fixed parameter of the EOSSP, thus reducing the potential target observations or data downlink through the restricted VTW access to targets and ground stations (when considered), respectively. The limitations can be alleviated using current state-of-the-art satellite concepts of operations (CONOPS), allowing the VTW to be manipulated as an intermediate decision variable \cite{LEMAITRE2002}. 

One such state-of-the-art satellite CONOPS is \textit{satellite agility}, a key aspect of the \textit{agile Earth observation satellite scheduling problem} (AEOSSP). Satellite agility refers to the ability of satellites to perform attitude control (slewing) in roll, pitch, and/or yaw \cite{LEMAITRE2002}. The AEOSSP leverages slewing as an additional task, allowing the satellites to change their pointing angle beyond the nadir direction, thus extending existing VTWs or creating new VTWs that were previously beyond the nadir directional swath width \cite{CHATTERJEE2022}. Prominent formulations of the AEOSSP include the use of slewing to maximize the number or profit of observations under various operational conditions \cite{Eddy2021,CHATTERJEE2022,Lee2024Optimal,Wang2019Scheduling} and the use of slewing to prioritize the amount of data returned to a ground station or operator \cite{Kim2020TaskScheduling,Herrmann2023Reinforcement,Herrmann2022}. However, satellite agility suffers from a degradation in observation quality due to the optical tilt and resolution change that results from angled observations, as a large slewing angle causes a difference in the look angle compared to nadir-directional observations \cite{Peng2020timedependent}. 

An additional state-of-the-art satellite CONOPS that has not yet been implemented as an extension to the EOSSP is \textit{constellation reconfigurability}. Constellation reconfigurability is defined as the capability of satellites within a constellation to perform orbital maneuvers to reform the constellation into a more optimal state \cite{deweck2008optimal,lee2023regional}. Such a capability provides flexibility and responsiveness to satellite systems \cite{paek2019,chen2015}, potentially allowing new opportunities for tasks that were previously unavailable. General uses of constellation reconfigurability external to EO include telecommunications systems with staged satellite deployment for the minimization of costs \cite{deweck2008optimal,Anderson2022megaconstellation,McGrath2015Design}, as well as response to incapacitated assets to compensate for performance losses \cite{Zuo2022Surrogate}. Investigations of constellation reconfigurability within EO include single-stage \cite{lee2023regional,Lee2021lagrangian,Lee2020binary} and multi-stage \cite{Lee2022,Lee2023,lee2024deterministic} reconfiguration for natural disaster impact monitoring, reconfiguration between selected operation modes (regional and global observation) \cite{paek2019,he2020reconfigurable}, and altitude changes both for observations of earthquake impacts \cite{McGrath2019General} and tracking mobile targets \cite{Morgan2023}. Within their investigations, Refs.~\cite{McGrath2019General,Morgan2023,McGrath2015Design} make use of low-thrust maneuvering, which, while valuable for fuel-efficient constellation reconfiguration, is computationally expensive and may not be as effective as high-thrust maneuvering in highly dynamic environments. Furthermore, previous research in Refs.~\cite{Pearl2023Comparing,Pearl2024Benchmarking} provides a comparative analysis between nadir-directional, agile, and maneuverable satellites within respective constellations relative to the obtained VTWs by each constellation of $100$ historical tropical cyclones. Within the comparative analysis put forth by Ref.~\cite{Pearl2024Benchmarking}, the best case of constellation reconfigurability outperforms the nadir-directional and agile satellite constellations by \SI{450.92}{\%} and \SI{48.52}{\%} on average, respectively, thus depicting promising results regarding the performance of constellation reconfigurability. Each application of constellation reconfigurability has demonstrated beneficial results through performance levels significantly greater than those obtained by fixed constellation configurations. Such research highlights the value that constellation reconfigurability provides to constellation systems; however, each implementation focuses on either maximum observations (maximum VTW), minimum revisit times, or a combination of both, hence neglecting data capacity, battery usage, and data downlink, all of which are key components of the EOSSP. 

This paper implements constellation reconfigurability to the EOSSP to develop the \textit{Reconfigurable Earth Observation Satellite Scheduling Problem} (REOSSP). Constellation reconfigurability is enabled through recent developments in the \textit{Multistage Constellation Reconfiguration Problem} (MCRP) \cite{Lee2022,Lee2023,lee2024deterministic} and is motivated through results indicating the performance level of constellation reconfigurability~\cite{Pearl2023Comparing,Pearl2024Benchmarking}. The novel REOSSP formulation is developed through the use of mixed-integer linear programming (MILP), enabling the acquisition of provably optimal solutions and the use of commercial optimization problem-solving software. Additionally, the presented formulation incorporates the consideration of onboard data storage, a feature neglected in Refs.~\cite{Wang2019Scheduling,Candela2023Dynamic,Swope2024Storm}, as well as onboard battery storage, a feature neglected in Ref.~\cite{Kim2020TaskScheduling}, whereas Refs.~\cite{Eddy2021,LEMAITRE2002,Chen2018MILP} neglect both features. Inspiration for the inclusion of such rigorous considerations originates from Refs.~\cite{Cho2018Optimization,Herrmann2023Reinforcement,Wang2023DeepRL} that consider these crucial aspects of the EOSSP.

Furthermore, this paper seeks to improve the computation runtime of the REOSSP through the use of a rolling horizon procedure (RHP). The RHP approach iteratively optimizes decisions over a given time horizon \cite{Morgan2014} by splitting a large problem into multiple subproblems containing a smaller portion of the full problem. As such, each subproblem has a reduced solution space and thus a reduced level of complexity. This paper presents the computational efficiency enabled by solving the REOSSP through the use of RHP.

In parallel, a novel EOSSP formulation, similarly considering onboard data and battery storage over time, is developed as a baseline nonreconfigurable equivalent to the REOSSP to benchmark the performance of the REOSSP through experimentation. Such experiments include a set of random instances with randomized parameters and a case study applied to Hurricane Sandy. Through these experiments and the results obtained, the conclusion of Ref.~\cite{Pearl2025Developing} that constellation reconfigurability greatly outperforms a standard nadir-directional constellation without reconfigurability is reinforced. Therefore, this paper additionally serves as further proof of concept for the implementation of constellation reconfigurability to the EOSSP at a higher level of fidelity and operational complexity than previous research. 

\begin{figure}[!ht]
    \centering
    \includegraphics[width=0.95\linewidth]{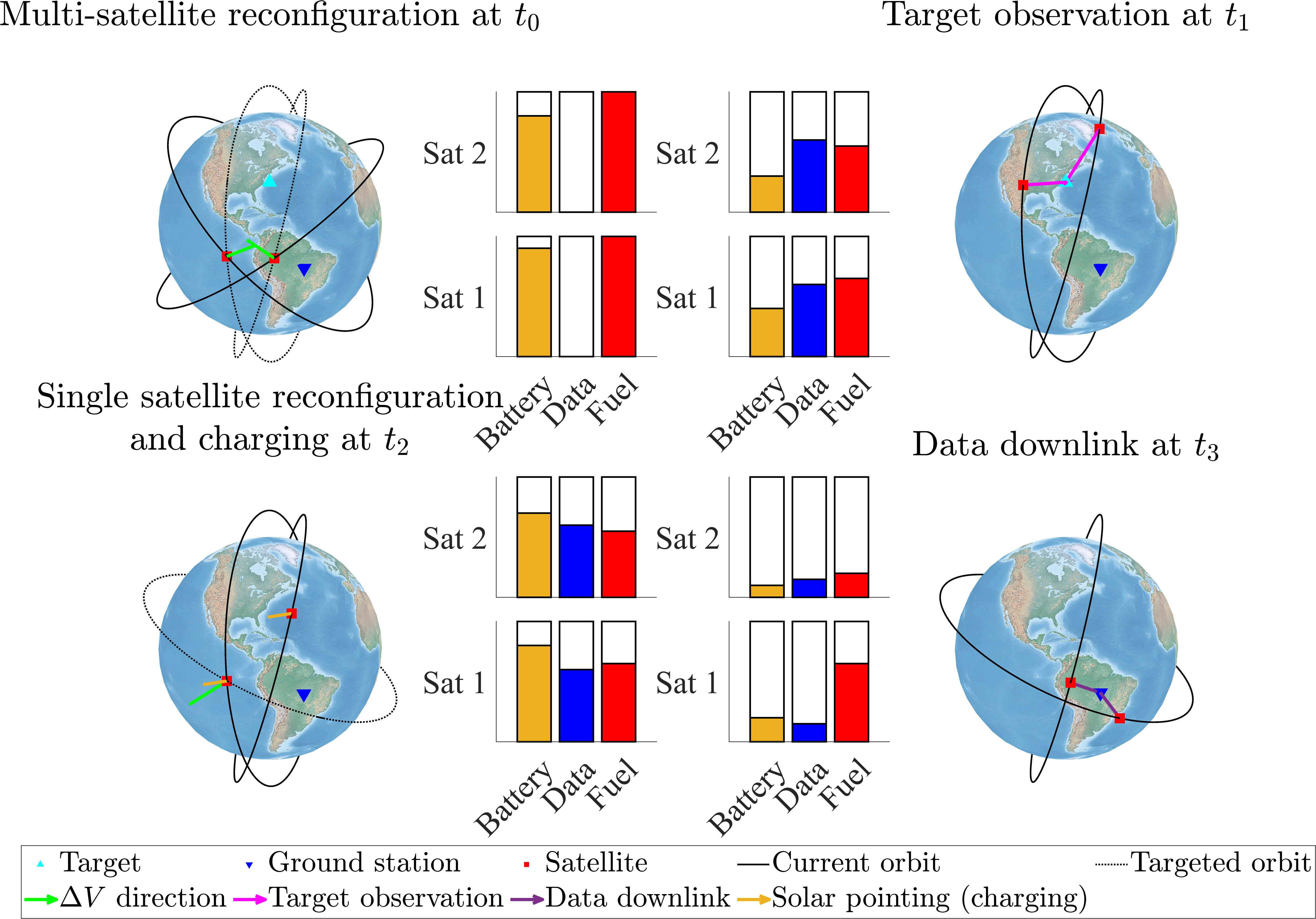}
    \caption{A demonstration of the capabilities enabled in the REOSSP.}
    \label{fig:REOSSP_demonstration}
\end{figure}

An illustrative example of the capabilities of the REOSSP is shown in Fig.~\ref{fig:REOSSP_demonstration}, including a representation of all tasks and operational constraints. Tasks, such as orbital maneuvering, target observation, data downlink to a ground station, and solar charging, are shown through arrows in the direction of the task, while the state of variables relating to operational constraints, such as the battery, data, and fuel capacities, are shown as a level that represents the current capacity. The initial state of the constellation is at $t_0$, where both satellites have no visibility of a given target. As a result of no present visibility, both satellites perform orbital maneuvers to form a more optimal configuration, a process that drains both the onboard battery and fuel. The reconfigured state of the constellation is at $t_1$, which has now gained visibility of the target, thus observing the target and contributing to the onboard data storage while further draining the onboard battery. At $t_2$, an observation of the target has been completed and may not be required to be performed again, leading both satellites to replenish the onboard battery storage through solar charging. Additionally, one satellite in this portion will not have visibility of the ground station, and as such, this satellite performs an additional orbital maneuver to gain said visibility. Finally, $t_3$ shows the final state of the constellation with both satellites downlinking the observed data to the ground station, transmitting the data currently stored at the cost of more battery power. Furthermore, this portion shows the additional drain on the fuel of the satellite that performed an additional orbital maneuver. 

The rest of this paper is organized as follows. Sections~\ref{sec:EOSSP} and~\ref{sec:REOSSP} introduce the MILP formulations of the \EOSSP and \REOSSP, respectively. Then, Sec.~\ref{sec:RHP} introduces an algorithmic solution method to the \REOSSP using RHP, a novel formulation referred to as the \RHP. Next, Sec.~\ref{sec:Experiments} conducts a thorough comparative analysis between the \EOSSP, \REOSSP, and \RHP, including instances with randomized parameters and a real-world case study using historical data from Hurricane Sandy. Finally, Sec.~\ref{sec:Conclusion} provides a review of overall results and contributes several possible directions for future research.

\section{Earth Observation Satellite Scheduling Problem---Baseline} \label{sec:EOSSP}

To incorporate constellation reconfigurability into the EOSSP, a baseline formulation is developed to be expanded for reconfiguration, which also proves useful as a comparison between a reconfigurable and a nonreconfigurable solution. This formulation is henceforth denoted as the \EOSSP. The formulation of the \EOSSP in this paper is an optimization problem making use of MILP to schedule tasks of target observation, data downlink to ground stations, and solar panel charging. Each task is constrained by physical capabilities, data storage availability, and/or battery power availability. 

The overall finite schedule duration is designated as $T_r$, which is discretized by a time step size $\Delta t>0$ to result in the finite number of discrete time steps $T=T_r/\Delta t$. The set of time steps, defined as $\mathcal{T} = \{1,2,\ldots,T\}$, contains all time steps at which a task can be performed. There also exists a set of satellites $\mathcal{K} = \{1,2,\ldots,K\}$ containing $K$ total satellites and the associated classical orbital elements, a set of targets for observation and associated positions defined as $\mathcal{P} = \{1,2,\ldots,P\}$ where $P$ is the total number of targets, and a set of ground stations for the downlink of data and associated positions defined as $\mathcal{G} = \{1,2,\ldots,G\}$ where $G$ is the total number of ground stations. 

The \EOSSP also utilizes various parameters within the problem constraints. The visibility of target $p \in \mathcal{P}$ by satellite $k \in \mathcal{K}$ at time step $t \in \mathcal{T}$ is contained in $V^k_{tp}$ as a binary value dependent on if $p$ is visible (one) or not (zero). Similarly, the visibility of ground station $g \in \mathcal{G}$ by satellite $k$ at time step $t$ is contained in $W^k_{tg}$ and the visibility of the sun by satellite $k$ at time step $t$ is contained in $H^k_t$, both of which are binary values.

Finally, the data and battery storage levels are tracked for each satellite $k$ from each time step $t$ to the next $t+1$. The current data storage level is determined by tracking the data gained through target observations and the data transmitted through downlink to any ground station. Additionally, the data storage level of each satellite $k$ is restricted to not deplete below a minimum threshold $D^k_{\min}$ and not accumulate more than the data storage capacity of the satellite, defined as $D^k_{\max}$. Similarly, the current battery storage level is determined by tracking the energy drained by observation, data downlink, and standard operational tasks such as tracking telemetry and time, and the energy replenished by charging through sunlight exposure. The battery storage level of each satellite $k$ is also restricted to not deplete below a minimum threshold $B^k_{\min}$ and not accumulate more than the battery capacity, defined as $B^k_{\max}$. 

\subsection{Decision Variables and Indicator Variables}

The decision variables of the \EOSSP are the tasks being scheduled throughout the duration of a given schedule, with additional indicator variables that track the data storage and battery levels over time. Each decision variable task is defined as binary such that a value of one corresponds to the task being performed and a value of zero corresponds to the task not being performed, while the indicator variables are defined as real numbers with upper and lower bounds defined as the maximum and minimum storage capacities, respectively. Each of these variables is applied to each satellite $k$ and time step $t$ in the overall schedule horizon, while the tasks additionally apply to their associated objective. 

The decision variables relative to satellite $k$ at time step $t$ include the observation of target $p$, $y^k_{tp}$ in constraints~\eqref{EOSSP:y}; the downlink of data to ground station $g$, $q^k_{tg}$ in constraints~\eqref{EOSSP:q}; and solar charging, $h^k_t$ in constraints~\eqref{EOSSP:h}. Similarly, the indicator variables relative to satellite $k$ at time step $t$ include: the current data storage level, $d^k_t$ in constraints~\eqref{EOSSP:D}, restricted between a minimum of $D^k_{\min}$ and a maximum of $D^k_{\max}$; and the current battery storage level, $b^k_t$ in constraints~\eqref{EOSSP:B}, restricted between a minimum of $B^k_{\min}$ and a maximum of $B^k_{\max}$.

\begin{subequations}
    \begin{alignat}{2}
        y^k_{tp} & \in \{0, 1\}, \quad && t \in \mathcal{T}, p \in \mathcal{P}, k \in \mathcal{K} \label{EOSSP:y} \\
        q^k_{tg} & \in \{0, 1\}, \quad && t \in \mathcal{T}, g \in \mathcal{G}, k \in \mathcal{K} \label{EOSSP:q} \\
        h^k_t & \in \{0, 1\}, \quad && t \in \mathcal{T}, k \in \mathcal{K} \label{EOSSP:h} \\
        d^k_t & \in [D^k_{\min}, D^k_{\max}], \quad && t \in \mathcal{T}, k \in \mathcal{K} \label{EOSSP:D} \\
        b^k_t & \in [B^k_{\min}, B^k_{\max}], \quad && t \in \mathcal{T}, k \in \mathcal{K} \label{EOSSP:B}
    \end{alignat}
    \label{EOSSP:dv}
\end{subequations}

\subsection{Objective Function}

The objective function of the \EOSSP maximizes the number of observations and subsequent downlink occurrences of data to ground stations, more heavily prioritizing the downlink of data through an arbitrary weight, $C > 1$. The objective function value of the \EOSSP is denoted as $z_{E}$ obtained via 
\begin{equation}
    z_{E} = \sum_{k \in \mathcal{K}} \sum_{t \in \mathcal{T}} \left( \sum_{g \in \mathcal{G}} C q^k_{tg} + \sum_{p \in \mathcal{P}} y^k_{tp} \right)
    \label{EOSSP:obj}
\end{equation}
where $z_{E}$ is determined through the weighted number of downlink tasks performed by the decision variable $q^k_{tg}$ and the unweighted number of observation tasks performed by the decision variable $y^k_{tp}$ with respect to all $K$ satellites, all $T$ time steps, and all $G$ ground stations or all $P$ targets for decision variables $q^k_{tg}$ and $y^k_{tp}$, respectively. Alternatively, the figure of merit for the \EOSSP corresponds to the amount of data downlinked to ground stations, obtained via 
\begin{equation}
    Z_{E} = \sum_{k \in \mathcal{K}} \sum_{t \in \mathcal{T}} \sum_{g \in \mathcal{G}} D_{\text{comm}} q^k_{tg} 
    \label{EOSSP:merit}
\end{equation}
where $D_{\text{comm}} \ge 0$ is the amount of data transmitted during a single time step. The objective function value and the figure of merit are both used to compare the performance of the \EOSSP to the \REOSSP and \RHP. 

\subsection{Constraints}

The \EOSSP contains constraints, including the application of target visibility, ground station visibility, sun visibility, maximum (and minimum) data storage capacity, and maximum (and minimum) battery storage capacity. 

\subsubsection{Time Window Constraints}

Each task is restricted to only be performed when visibility of the task is available to each satellite individually; these constraints define the VTW for each task utilizing the associated parameter such as $V^k_{tp}$ for targets observation, $W^k_{tg}$ for ground station downlink, and $H^k_t$ for solar panel charging, respectively. Additionally, power-related tasks such as target observation, data downlink, and solar charging are restricted to occurring exclusively from one another.
\begin{subequations}
    \begin{alignat}{2}
& V^k_{tp} \ge y^k_{tp}, \quad && \forall t \in \mathcal{T}, \forall p \in \mathcal{P}, \forall k \in \mathcal{K}
\label{EOSSP:target_visibility}\\
& W^k_{tg} \ge q^k_{tg}, \quad && \forall t \in \mathcal{T}, \forall g \in \mathcal{G}, \forall k \in \mathcal{K}
\label{EOSSP:gs_visibility}\\
& H^k_t \ge h^k_t, \quad && \forall t \in \mathcal{T}, \forall k \in \mathcal{K}
\label{EOSSP:sun_visibility}\\
& \sum_{p \in \mathcal{P}} y^k_{tp} + \sum_{g \in \mathcal{G}} q^k_{tg} + h^k_t \leq 1, \quad && \forall t \in \mathcal{T}, \forall k \in \mathcal{K}
\label{EOSSP:obvs-down_overlap}
    \end{alignat}
    \label{EOSSP:Visibility}
\end{subequations}

Constraints~(\ref{EOSSP:target_visibility}--\ref{EOSSP:sun_visibility}) restrict the observation of targets, downlink of data to ground stations, and solar charging to the associated VTW, respectively. These constraints operate under the assumption that only one satellite is required to perform each task such that each satellite operates independently. Additionally, constraints~\eqref{EOSSP:obvs-down_overlap} restrict observations, data downlink, and solar charging only to be performed in exclusion to one another, ensuring that each satellite must choose which task to perform if the VTW of each task overlaps. Such a consideration is in place for battery routing consideration to not overload the simultaneous power exchange within an onboard battery, similar to the task overlap exclusion employed in Ref.~\cite{Huang2021}. Additionally, these constraints ensure that only one target may be viewed at each time step or that only one ground station may be utilized for data downlink for the same reason. 

\begin{figure}[!ht]
    \centering
    \includegraphics[width = 0.5\textwidth]{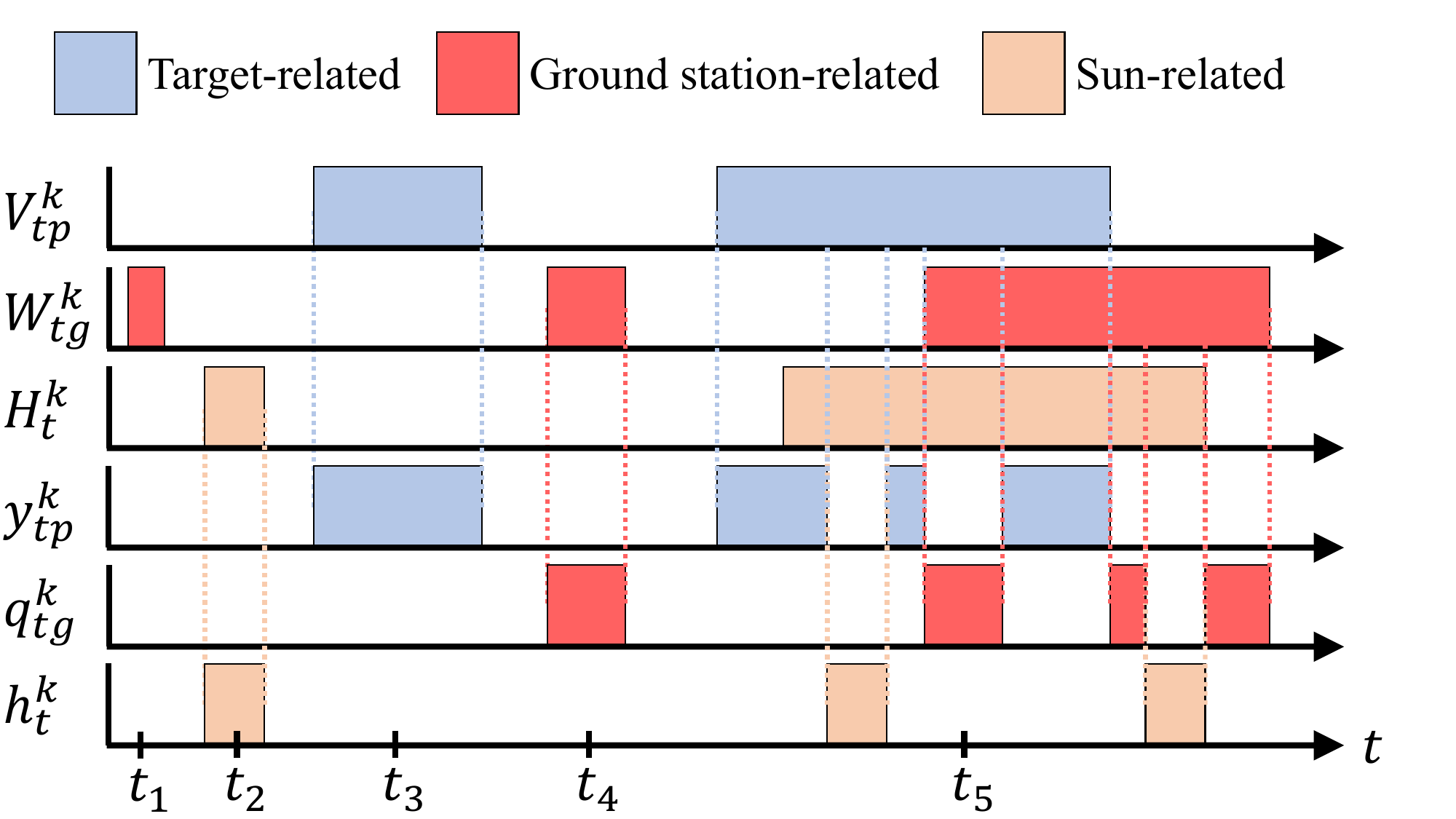}
    \caption{Feasible observation, downlink, and solar charging occurrence.}
    \label{fig:EOSSP_y_vs_q}
\end{figure}

A demonstration of the visibility conditions in constraints~(\ref{EOSSP:target_visibility}--\ref{EOSSP:sun_visibility}) as well as the overlap exclusion conditions in constraints~\eqref{EOSSP:obvs-down_overlap} is shown in Fig.~\ref{fig:EOSSP_y_vs_q}. Within the figure, illustrative binary visibility parameters $V^k_{tp}$, $W^k_{tg}$, and $H^k_t$ are shown through blue, red, and orange boxes, respectively, indicating time along a horizon with visibility of the associated goal, and the binary decision variables $y^k_{tp}$, $q^k_{tg}$, and $h^k_t$ are shown similarly. The figure depicts some key occurrences resulting from the above constraints, labeled as time steps $t_1$ through $t_5$. First, time step $t_1$ depicts the visibility of a ground station with no downlink occurrence as a result of no prior observations contributing data to be downlinked. Next time step $t_2$ depicts the visibility of the sun along with the subsequent solar charging. Then, time step $t_3$ shows the visibility and subsequent observation of a target, and time step $t_4$ illustrates the visibility of a ground station that allows the downlink of the previously observed data. Finally, time step $t_5$ depicts the enforcement of constraints~\eqref{EOSSP:obvs-down_overlap}, restricting only one task to occur at each time step, despite all three tasks having visibility of the associated goal at various times. The enforcement of constraints~\eqref{EOSSP:obvs-down_overlap} additionally applies to the simultaneous visibility of different targets and ground stations, though this is not shown in the figure directly. 

\subsubsection{Data Tracking and Storage Constraints}

The second set of constraints tracks the usage of each satellite's onboard data storage, determining when data is gained or downlinked and restricting the data so that it does not exceed storage limits.
\begin{subequations}
    \begin{alignat}{2}
& d^k_{t+1} = d^k_t + \sum_{p \in \mathcal{P}} D_{\text{obs}} y^k_{tp} - \sum_{g \in \mathcal{G}} D_{\text{comm}} q^k_{tg} , \quad && \forall t \in \mathcal{T} \setminus \{T\}, \forall k \in \mathcal{K}
\label{EOSSP:d-track}\\
& d^k_t + \sum_{p \in \mathcal{P}} D_{\text{obs}} y^k_{tp} \leq D^k_{\max}, \quad && \forall t \in \mathcal{T}, \forall k \in \mathcal{K}
\label{EOSSP:d<max}\\
& d^k_t - \sum_{g \in \mathcal{G}} D_{\text{comm}} q^k_{tg} \ge D^k_{\min}, \quad && \forall t \in \mathcal{T}, \forall k \in \mathcal{K}
\label{EOSSP:d>0}
    \end{alignat}
    \label{EOSSP:Data}
\end{subequations}

Constraints~\eqref{EOSSP:d-track} track the data storage of each satellite at each time step, ensuring that the tasks performed at time step $t$ contribute to the data storage level at time step $t+1$, where $D_{\text{obs}} \ge 0$ is the amount of data acquired through one observation within one time step. Furthermore, $D_{\text{obs}}$ is a constant value of given units, while $D_{\text{comm}}$ is obtained by multiplying a data transmission rate in units/time by $\Delta t$ to obtain a discrete constant value. It is also assumed that each satellite begins with the minimum data storage value such that $d^k_1 = D^k_{\min},~\forall k \in \mathcal{K}$, allowing the first data observed to contribute at time $t=1$ for the data storage level of $d^k_2$. Constraints~\eqref{EOSSP:d<max} ensure that the data stored within each satellite through observation does not exceed the maximum capacity set by $D^k_{\max}$. Similarly, constraints~\eqref{EOSSP:d>0} ensure that the data transmitted to a ground station does not exceed the actual amount of data contained by the satellite. 

A demonstration of the control and contribution of the data storage for each satellite is shown in Fig.~\ref{fig:EOSSP_data}, with similar key occurrences denoted by time steps $t_1$ through $t_5$. First, time steps $t_1$ and $t_2$ depict the contribution of observations and the drain of downlink to the onboard data storage, respectively, provided through constraints~\eqref{EOSSP:d-track}. Then, time steps $t_3$ and $t_4$ depict the halt of observations such that the onboard data level does not exceed the maximum capacity set by $D^k_{\max}$, provided through constraints~\eqref{EOSSP:d<max}. Finally, time step $t_5$ depicts the halt of downlink once the data level is fully drained, ensuring the current data level remains non-negative, provided through constraints~\eqref{EOSSP:d>0}.

\subsubsection{Battery Tracking and Storage Constraints}

The third and final set of constraints tracks the usage of the onboard battery power of each satellite, similarly to that of the data tracking constraints, determining when the battery is charged or depleted, as well as restricting the power to not overflow onboard battery limits and not deplete below the minimum capacity.
\begin{subequations}
    \begin{alignat}{2}
& b^k_{t+1} = b^k_t + B_{\text{charge}} h^k_t - \sum_{p \in \mathcal{P}} B_{\text{obs}} y^k_{tp} - \sum_{g \in \mathcal{G}} B_{\text{comm}} q^k_{tg} - B_{\text{time}}, \quad && \forall t \in \mathcal{T} \setminus \{T\}, \forall k \in \mathcal{K}
\label{EOSSP:b-track}\\
& b^k_t + B_{\text{charge}} h^k_t \leq B^k_{\max}, \quad && \forall t \in \mathcal{T}, \forall k \in \mathcal{K}
\label{EOSSP:b<max}\\
& b^k_t - \sum_{p \in \mathcal{P}} B_{\text{obs}} y^k_{tp} - \sum_{g \in \mathcal{G}} B_{\text{comm}} q^k_{tg} - B_{\text{time}} \ge B^k_{\min}, \quad && \forall t \in \mathcal{T}, \forall k \in \mathcal{K}
\label{EOSSP:b>0}
    \end{alignat}
    \label{EOSSP:Battery}
\end{subequations}

Constraints~\eqref{EOSSP:b-track} track the battery level of each satellite at each time step, where $B_{\text{charge}} \ge 0$ is the amount of power able to be charged, $B_{\text{obs}} \ge 0$ is the amount of power required for observation to occur, $B_{\text{comm}} \ge 0$ is the amount of power required for a data downlink transmission, and $B_{\text{time}} \ge 0$ is the power required for normal operational tasks such as telemetry and timekeeping. Each power value accounts for one time step and is found by multiplying a charge rate in units/time by $\Delta t$ to obtain a power value in the given units. It is similarly assumed that each satellite begins with a full battery charge such that $b^k_1 = B^k_{\max},~\forall k \in \mathcal{K}$, allowing the first power draw occurring at $t=1$ to contribute to the battery storage level of $b^k_2$. Constraints~\eqref{EOSSP:b<max} ensure that the power contained within each satellite through charging does not exceed the maximum capacity of the battery set by $B^k_{\max}$. Similarly, constraints~\eqref{EOSSP:b>0} ensure that the power required by other task performance does not drain the battery below the minimum capacity set by $B^k_{\min}$. 

\begin{figure}[!ht]
    \centering
    \begin{subfigure}[h]{0.49\textwidth}
        \centering
        \includegraphics[width = \textwidth]{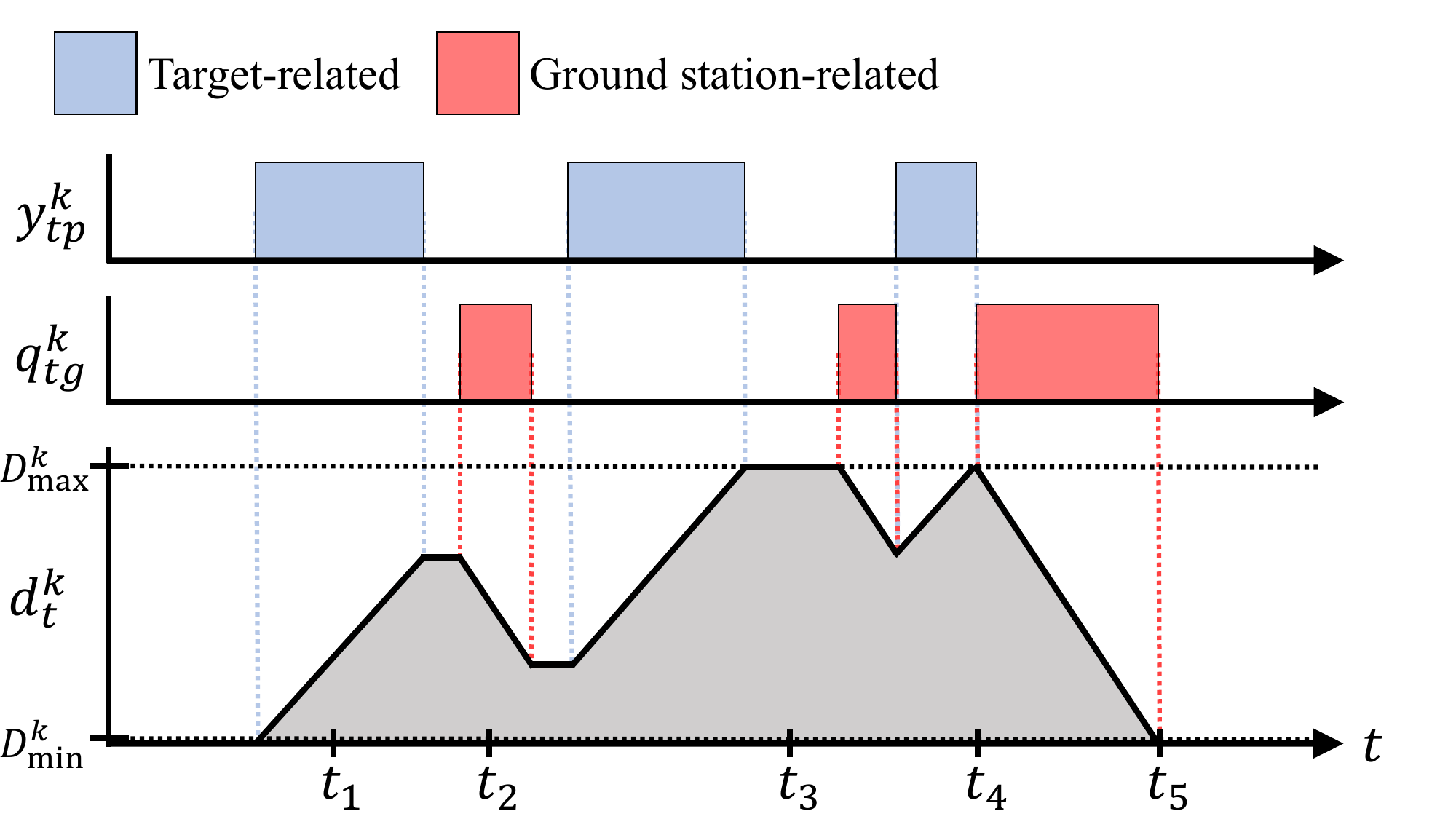}
        \caption{Visualization of the data storage level of satellite $k$}
        \label{fig:EOSSP_data}
    \end{subfigure}
    \hfill
    \begin{subfigure}[h]{0.49\textwidth}
        \centering
        \includegraphics[width = \textwidth]{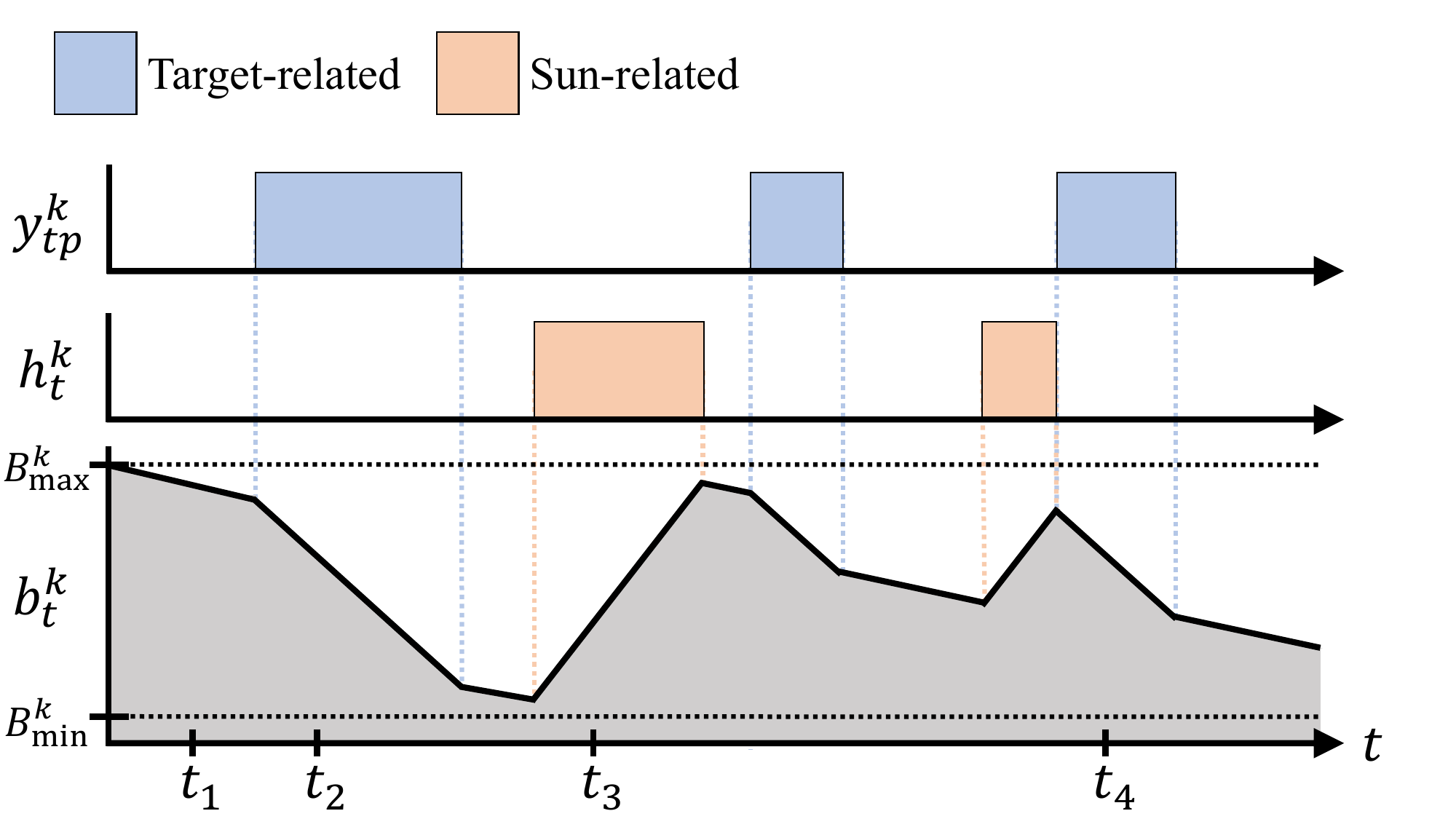}
        \caption{Visualization of the battery storage level of satellite $k$}
        \label{fig:EOSSP_battery}
    \end{subfigure}
    \caption{Data and battery storage level visualizations of satellite $k$.}
\end{figure}

Figure~\ref{fig:EOSSP_battery} visualizes the battery consumption similarly to Fig.~\ref{fig:EOSSP_data}, with key occurrences denoted by time steps $t_1$ through $t_4$. First, time step $t_1$ demonstrates the constant battery drain caused by normal operational tasks, $B_{\text{time}}$, and time step $t_2$ demonstrates the combined battery drain of the constant drain and observation. Then, time step $t_3$ demonstrates the battery charging as a result of charging operations while within sun visibility. Finally, time step $t_4$ reiterates the task overlap exclusion from constraints~\eqref{EOSSP:obvs-down_overlap}. Such consumption of battery power additionally applies to data downlink, though it is not shown directly in the figure. 

\subsection{Full Formulation}

Given the objective function, constraints, and decision/indicator variables, the full \EOSSP is defined as follows: \hypertarget{EOSSP}{}
\begin{equation}
    \begin{split}
        \max & \quad z_{E} \\
        \text{s.t.} & \quad \text{Constraints \eqref{EOSSP:dv}--\eqref{EOSSP:obj}, \eqref{EOSSP:Visibility}--\eqref{EOSSP:Battery}}
    \end{split} \tag{\textsf{EOSSP}}
    \label{EOSSP:EOSSP}
\end{equation}

The objective function balances the number of occurrences of data downlink to any ground stations $g \in \mathcal{G}$ with the less heavily weighted number of occurrences of the observation of any targets $p \in \mathcal{P}$, while the constraints restrict the problem to the observation of data before data downlink, the assurance of available data storage, and the assurance of available power for task completion. Additionally, the \EOSSP provides the schedule of charging, observation, and data downlink performed by a constellation of satellites over a given time horizon. 

\section{Reconfigurable Earth Observation Satellite Scheduling Problem} \label{sec:REOSSP}

The REOSSP formulation, henceforth referred to as the \REOSSP, is an extension of the \EOSSP that incorporates satellite maneuvering tasks in addition to the tasks present in the \EOSSP, thus allowing a new degree of freedom. The constellation reconfiguration process is derived from the MCRP, obtained from Refs.~\cite{Lee2022,Lee2023,lee2024deterministic}. All sets and parameters remain as defined in Sec.~\ref{sec:EOSSP} unless otherwise specified.

Initial changes are made to $\mathcal{T}$ regarding stages for reconfiguration, as well as additional orbital slots being provided to each satellite $k \in \mathcal{K}$. Firstly, within $\mathcal{T}$ there exists a set of equally spaced stages at which each satellite in the constellation can perform orbital maneuvers, defined as $\mathcal{S} = \{0,1,2,\ldots,S\}$, where $S$ is the number of equally spaced stages and stage $s=0$ is the stage wherein the initial configuration of the constellation is present. As a result of each stage being equally spaced in $\mathcal{T}$, there exists a set of time steps within each stage defined as $\mathcal{T}^s = \{1,2,\ldots,T^s\}$ where $T^s=T/S$ is an integer defining the number of time steps per stage. Secondly, the orbital slots are defined as $\mathcal{J}^{sk} = \{1,2,\ldots,J^{sk}\}$, where $J^{sk}$ is the total number of available slots for satellite $k \in \mathcal{K}$ and stage $s \in \mathcal{S}$. 

As a result of the \REOSSP dividing the time horizon into stages, some parameters are modified slightly, and additional parameters are included for the cost of orbital maneuvers. Most importantly, the binary visibility parameters for targets, ground stations, and the sun are extended to include two additional dimensions representing the current stage and occupied orbital slot; each is now defined as $V^{sk}_{tjp}$, $W^{sk}_{tjg}$, and $H^{sk}_{tj}$, respectively. Additionally, the cost of transfer between two orbital slots of a given satellite and a given stage is $c^{sk}_{ij}$, defined as the transfer between slot $i \in \mathcal{J}^{s-1, k}$ in the previous stage configuration ($s-1$) to slot $j \in \mathcal{J}^{sk}$ in the current stage ($s$) relative to satellite $k$ and stage $s \in \mathcal{S} \setminus \{0\}$. Furthermore, the maximum transfer cost for the entire schedule is allotted to each satellite as $c^k_{\max},~\forall k \in \mathcal{K}$.

\subsection{Decision Variables}

The decision variables and indicator variables of the \REOSSP are the same tasks and indicators included in the \EOSSP, with additional dimensionality regarding reconfiguration stages and an additional task to control the maneuvers of each satellite. The new decision variable, $x^{sk}_{ij}$ in constraints~\eqref{REOSSP:x}, defines the sequence of maneuvers of satellite $k$ from an orbital slot $i \in \mathcal{J}^{s-1, k}$ in the previous stage to an orbital slot $j \in \mathcal{J}^{sk}$ in the current stage for all stages $s \in \mathcal{S} \setminus \{0\}$ such that the singleton set initial condition orbital slot option, $\mathcal{J}^{0k}$, is transferred from stage $s=1$. All other decision and indicator variables expand the dimensionality of the respective \EOSSP variables through an additional index $s \in \mathcal{S} \setminus \{0\}$ relative to time step $t \in \mathcal{T}^s$ and are represented in constraints~(\ref{REOSSP:y}--\ref{REOSSP:D}):
\begin{subequations}
    \begin{alignat}{2}
        x^{sk}_{ij} & \in \{0, 1\}, \quad && s \in \mathcal{S} \setminus \{0\}, k \in \mathcal{K}, i \in \mathcal{J}^{s-1, k}, j \in \mathcal{J}^{sk} \label{REOSSP:x} \\
        y^{sk}_{tp} & \in \{0, 1\}, \quad && s \in \mathcal{S} \setminus \{0\}, t \in \mathcal{T}^s, p \in \mathcal{P}, k \in \mathcal{K} \label{REOSSP:y} \\
        q^{sk}_{tg} & \in \{0, 1\}, \quad && s \in \mathcal{S} \setminus \{0\}, t \in \mathcal{T}^s, g \in \mathcal{G}, k \in \mathcal{K} \label{REOSSP:q} \\
        h^{sk}_t & \in \{0, 1\}, \quad && s \in \mathcal{S} \setminus \{0\}, t \in \mathcal{T}^s, k \in \mathcal{K} \label{REOSSP:h} \\
        d^{sk}_t & \in [D^k_{\min}, D^k_{\max}], \quad && s \in \mathcal{S} \setminus \{0\}, t \in \mathcal{T}^s, k \in \mathcal{K} \label{REOSSP:D} \\
        b^{sk}_t & \in [B^k_{\min}, B^k_{\max}], \quad && s \in \mathcal{S} \setminus \{0\}, t \in \mathcal{T}^s, k \in \mathcal{K} \label{REOSSP:B} 
    \end{alignat}
    \label{REOSSP:dv}
\end{subequations}

\subsection{Objective Function}

The objective function of the \REOSSP mirrors the \EOSSP, maximizing the number of observations and subsequent downlink occurrences of data to ground stations with the same weighting $C$ to prioritize the downlink of data. The objective function value of the \REOSSP is denoted as $z_{R}$ obtained via 
\begin{equation}
    z_{R} = \sum_{k \in \mathcal{K}} \sum_{s \in \mathcal{S} \setminus \{0\}} \sum_{t \in \mathcal{T}^s} \left( \sum_{g \in \mathcal{G}} C q^{sk}_{tg} + \sum_{p \in \mathcal{P}} y^{sk}_{tp} \right)
    \label{REOSSP:obj}
\end{equation}
where $z_{R}$ additionally considers all $S$ stages and all $T^s$ time steps per stage. Similarly to the \EOSSP, the figure of merit for the \REOSSP corresponds to the amount of data downlinked to ground stations, obtained via
\begin{equation}
    Z_{R} = \sum_{k \in \mathcal{K}} \sum_{s \in \mathcal{S} \setminus \{0\}} \sum_{t \in \mathcal{T}^s} \sum_{g \in \mathcal{G}} D_{\text{comm}} q^{sk}_{tg}
    \label{REOSSP:merit}
\end{equation}
The objective function value and the figure of merit are used to compare the performance of the \REOSSP to the \EOSSP and \RHP. To distinguish the resultant objective function values and figure of merit values, the subscript $R$ is used in place of the subscript $E$. 

\subsection{Constraints}

The \REOSSP contains similar constraints to the \EOSSP but with extended dimensionality to incorporate stages sequentially within the time horizon and reflect the effects of the reconfiguration process through new satellite positioning and resources required to perform such maneuvers. 

\subsubsection{Orbital Maneuver Path Continuity Constraints}

The first set of constraints restricts the path of constellation reconfigurability to feasible maneuvering sequences concerning a given budget for propellant.
\begin{subequations}
    \begin{alignat}{2}
& \sum_{j \in \mathcal{J}^{1k}} x^{1k}_{ij} = 1, \quad && \forall k \in \mathcal{K}, \forall i \in \mathcal{J}^{0k} 
\label{REOSSP:flow_from_initial_conditions}\\
& \sum_{j \in \mathcal{J}^{s+1, k}} x^{s+1, k}_{ij} - \sum_{j' \in \mathcal{J}^{s-1, k}} x^{sk}_{j'i} = 0, \quad && \forall s \in \mathcal{S} \setminus \{0, S\}, \forall k \in \mathcal{K}, \forall i \in \mathcal{J}^{sk} 
\label{REOSSP:flow_from0<s<S}\\
& \sum_{s \in \mathcal{S} \setminus \{0\}} \sum_{i \in \mathcal{J}^{s-1, k}} \sum_{j \in \mathcal{J}^{sk}} c^{sk}_{ij} x^{sk}_{ij} \leq c^k_{\max}, \quad && \forall k \in \mathcal{K} 
\label{REOSSP:flow_cost}
    \end{alignat}
    \label{REOSSP:flow_constraints}
\end{subequations}

Constraints~\eqref{REOSSP:flow_from_initial_conditions} ensure that only one orbital slot in the first stage ($\mathcal{J}^{1k}$) can be selected for transfer from the initial condition ($\mathcal{J}^{0k}$), ensuring that one satellite cannot occupy multiple positions at one time. Constraints~\eqref{REOSSP:flow_from0<s<S} similarly ensure that only one orbital slot in subsequent stages is selected for transfer but also extend the condition to state that the satellite can only transfer from the current slot ($i \in \mathcal{J}^{sk}$) if it arrived there previously. Both constraints ensure that satellites do not conduct transfers from positions in which they may not be located. Finally, constraints~\eqref{REOSSP:flow_cost} restrict the allowable transfer between slots to within the maximum propellant budget of the satellite.

\begin{figure}[!ht]
    \centering
    \includegraphics[width = 0.6\textwidth]{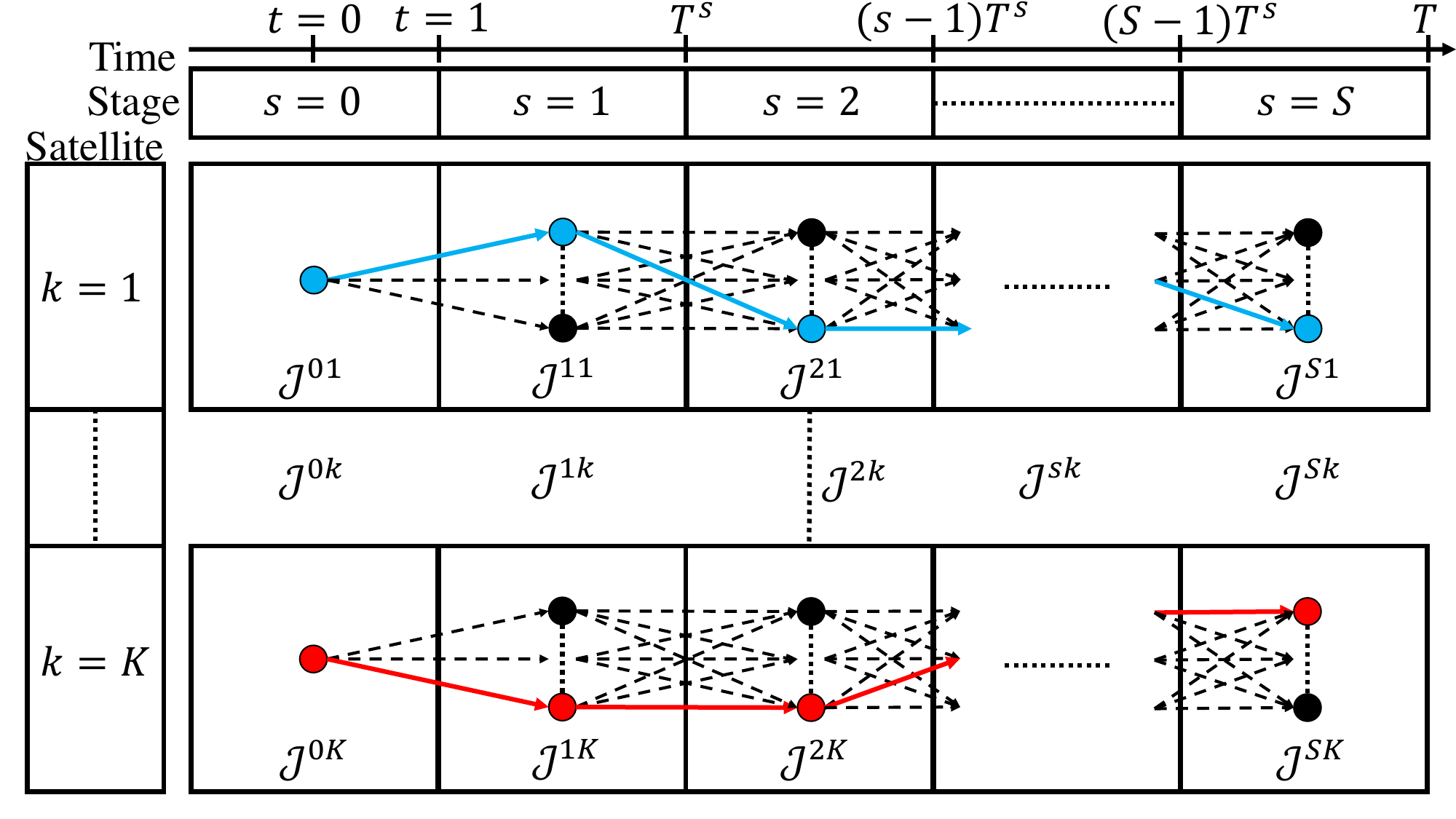}
    \caption{Feasible orbital maneuver progression for $K$ satellites.}
    \label{fig:recon_demonstration}
\end{figure}

A demonstration of feasible path continuity and maneuvering is shown in Fig.~\ref{fig:recon_demonstration}, depicting the possible paths taken by any number of satellites for any number of orbital slots or stages, where the highlighted paths represent $x^{sk}_{ij} = 1$. The figure notes that the initial condition and first orbital maneuvers occur outside of the time horizon such that the constellation configuration at stage $s=1$ is formed at time $t=1$ and subsequent stages occur at the corresponding $t=1+(s-1)T^s$ time. Additionally, the figure demonstrates the option to remain within the current orbital slot, such as satellite $K$ in stage two, resulting in a cost $c^{sk}_{ij}$ of zero. 

\subsubsection{Time Window Constraints}

The second set of constraints restricts task performance to the available VTWs, additionally enforcing the task overlap exclusion of each task. These constraints mirror constraints~\eqref{EOSSP:Visibility} while additionally accounting for the orbital maneuver path and subsequent satellite orbital slot positions within each reconfiguration stage, incorporated using $x^{sk}_{ij}$ and other extended dimensionality parameters.
\begin{subequations}
    \begin{alignat}{2}
& \sum_{i \in \mathcal{J}^{s-1, k}} \sum_{j \in \mathcal{J}^{sk}} V^{sk}_{tjp} x^{sk}_{ij} \ge y^{sk}_{tp}, \quad && \forall s \in \mathcal{S} \setminus \{0\}, \forall t \in \mathcal{T}^s, \forall p \in \mathcal{P}, \forall k \in \mathcal{K}
\label{REOSSP:target_visibility}\\
& \sum_{i \in \mathcal{J}^{s-1, k}} \sum_{j \in \mathcal{J}^{sk}} W^{sk}_{tjg} x^{sk}_{ij} \ge q^{sk}_{tg}, \quad && \forall s \in \mathcal{S} \setminus \{0\}, \forall t \in \mathcal{T}^s, \forall g \in \mathcal{G}, \forall k \in \mathcal{K}
\label{REOSSP:gs_visibility}\\
& \sum_{i \in \mathcal{J}^{s-1, k}} \sum_{j \in \mathcal{J}^{sk}} H^{sk}_{tj} x^{sk}_{ij} \ge h^{sk}_t, \quad && \forall s \in \mathcal{S} \setminus \{0\}, \forall t \in \mathcal{T}^s, \forall k \in \mathcal{K}
\label{REOSSP:sun_visibility}\\
& \sum_{p \in \mathcal{P}} y^{sk}_{tp} + \sum_{g \in \mathcal{G}} q^{sk}_{tg} + h^{sk}_t \leq 1, \quad && \forall s \in \mathcal{S} \setminus \{0\}, \forall t \in \mathcal{T}^s, \forall k \in \mathcal{K}
\label{REOSSP:obvs-down_overlap}
    \end{alignat}
    \label{REOSSP:Visibility}
\end{subequations}

\subsubsection{Data Tracking and Storage Constraints}

The third set of constraints tracks the usage of the onboard data storage level of each satellite, mirroring constraints~\eqref{EOSSP:Data}, while additionally splitting constraints~\eqref{EOSSP:d-track} between constraints~\eqref{REOSSP:d-track_not_Ts} and~\eqref{REOSSP:d-track_Ts}. Constraints~\eqref{REOSSP:d-track_not_Ts} track the data storage level of each satellite within the time horizon of each stage, while constraints~\eqref{REOSSP:d-track_Ts} track the data storage level at the stage gap, where the end of one stage meets the beginning of the next, providing the condition that $d^{sk}_{T^s+1} = d^{s+1, k}_1$.

\begin{subequations}
    \begin{alignat}{2}
& d^{sk}_{t+1} = d^{sk}_t + \sum_{p \in \mathcal{P}} D_{\text{obs}} y^{sk}_{tp} - \sum_{g \in \mathcal{G}} D_{\text{comm}} q^{sk}_{tg}, \quad && \forall s \in \mathcal{S} \setminus \{0\}, \forall t \in \mathcal{T}^s \setminus \{T^s\}, \forall k \in \mathcal{K}
\label{REOSSP:d-track_not_Ts}\\
& d^{s+1, k}_1 = d^{sk}_{T^s} + \sum_{p \in \mathcal{P}} D_{\text{obs}} y^{sk}_{T^s p} - \sum_{g \in \mathcal{G}} D_{\text{comm}} q^{sk}_{T^s g}, \quad && \forall s \in \mathcal{S} \setminus \{0, S\}, \forall k \in \mathcal{K}
\label{REOSSP:d-track_Ts}\\
& d^{sk}_t + \sum_{p \in \mathcal{P}} D_{\text{obs}} y^{sk}_{tp} \leq D^k_{\max}, \quad && \forall s \in \mathcal{S} \setminus \{0\}, \forall t \in \mathcal{T}^s , \forall k \in \mathcal{K}
\label{REOSSP:d<max}\\
& d^{sk}_t - \sum_{g \in \mathcal{G}} D_{\text{comm}} q^{sk}_{tg} \ge D^k_{\min}, \quad && \forall s \in \mathcal{S} \setminus \{0\}, \forall t \in \mathcal{T}^s , \forall k \in \mathcal{K}
\label{REOSSP:d>0}
    \end{alignat}
    \label{REOSSP:Data}
\end{subequations}

\subsubsection{Battery Tracking Constraints}

The fourth set of constraints only \textit{tracks} the onboard battery capacity, mirroring constraints~\eqref{EOSSP:b-track}. Similarly to constraints~\eqref{REOSSP:d-track_not_Ts} and~\eqref{REOSSP:d-track_Ts}, constraints~\eqref{REOSSP:b-track_not_Ts} and~\eqref{REOSSP:b-track_Ts_not_s1} track the battery storage level of each satellite within the time horizon of each stage and at the stage gap, respectively, where $B_{\text{recon}} \ge 0$ is the power cost of orbital maneuvers as a conservative estimate for thruster pointing. Separately, constraints~\eqref{REOSSP:b-track_Ts_s1} account for the first orbital maneuvers, defined as a reduction from a full onboard battery of $B^k_{\max}$.
\begin{subequations}
    \begin{alignat}{2}
& b^{sk}_{t+1} = b^{sk}_t + B_{\text{charge}} h^{sk}_t - \sum_{p \in \mathcal{P}} B_{\text{obs}} y^{sk}_{tp} - \sum_{g \in \mathcal{G}} B_{\text{comm}} q^{sk}_{tg} - B_{\text{time}}, \quad && \forall s \in \mathcal{S} \setminus \{0\}, \forall t \in \mathcal{T}^s \setminus \{T^s\}, \forall k \in \mathcal{K} 
\label{REOSSP:b-track_not_Ts}\\
\begin{split}
    &b^{s+1, k}_1 = b^{sk}_{T^s} + B_{\text{charge}} h^{sk}_{T^s} - \sum_{p \in \mathcal{P}} B_{\text{obs}} y^{sk}_{T^s p} - \sum_{g \in \mathcal{G}} B_{\text{comm}} q^{sk}_{T^s g} - \\
    & \qquad \sum_{i \in \mathcal{J}^{sk}} \sum_{j \in \mathcal{J}^{s+1, k}} B_{\text{recon}} x^{s+1, k}_{ij}  - B_{\text{time}}, 
\end{split} \quad && \forall s \in \mathcal{S} \setminus \{0, S\}, \forall k \in \mathcal{K}
\label{REOSSP:b-track_Ts_not_s1}\\
& b^{1k}_1 = B^k_{\max} - \sum_{i \in \mathcal{J}^{0k}} \sum_{j \in \mathcal{J}^{1k}} B_{\text{recon}} x^{1k}_{ij} , \quad && \forall k \in \mathcal{K}
\label{REOSSP:b-track_Ts_s1}
    \end{alignat}
    \label{REOSSP:Battery_Track}
\end{subequations}

\subsubsection{Battery Storage Constraints}

The final set of constraints ensures the onboard battery level remains within established bounds, where  constraints~\eqref{REOSSP:b<max} and~\eqref{REOSSP:b>0_not_Ts} are direct parallels to constraints~\eqref{EOSSP:b<max} and~\eqref{EOSSP:b>0}, respectively. Additionally, constraints~\eqref{REOSSP:b>0_Ts_not_s1} enforce the same condition as constraints~\eqref{REOSSP:b>0_not_Ts} at the stage gaps, and constraints~\eqref{REOSSP:b>0_Ts_s1} account for the first orbital maneuvers.

\begin{subequations}
    \begin{alignat}{2}
& b^{sk}_t + B_{\text{charge}} h^{sk}_t \leq B^k_{\max}, \quad && \forall s \in \mathcal{S} \setminus \{0\}, \forall t \in \mathcal{T}^s, \forall k \in \mathcal{K}
\label{REOSSP:b<max}\\
& b^{sk}_t - \sum_{p \in \mathcal{P}} B_{\text{obs}}\, y^{sk}_{tp} - \sum_{g \in \mathcal{G}} B_{\text{comm}}\, q^{sk}_{tg} - B_{\text{time}} \ge B^k_{\min}, \quad && \forall s \in \mathcal{S} \setminus \{0\}, \forall t \in 
\begin{cases}
    \mathcal{T}^s \setminus \{T^s\} & \text{if } s \neq S, \\
    \mathcal{T}^S & \text{if } s = S,
\end{cases} \quad \forall k \in \mathcal{K}
\label{REOSSP:b>0_not_Ts}\\
\begin{split}
& b^{sk}_{T^s} - \sum_{p \in \mathcal{P}} B_{\text{obs}} y^{sk}_{T^s p} - \sum_{g \in \mathcal{G}} B_{\text{comm}} q^{sk}_{T^s g} - \\
& \qquad \sum_{i \in \mathcal{J}^{sk}} \sum_{j \in \mathcal{J}^{s+1, k}} B_{\text{recon}} x^{s+1, k}_{ij} - B_{\text{time}} \ge B^k_{\min}, 
\end{split} \quad && \forall s \in \mathcal{S} \setminus \{0,S\}, \forall k \in \mathcal{K}
\label{REOSSP:b>0_Ts_not_s1}\\
& B^k_{\max} - \sum_{i \in \mathcal{J}^{0k}} \sum_{j \in \mathcal{J}^{1k}} B_{\text{recon}} x^{1k}_{ij} \ge B^k_{\min}, \quad && \forall k \in \mathcal{K}
\label{REOSSP:b>0_Ts_s1}
    \end{alignat}
    \label{REOSSP:Battery_Constrain}
\end{subequations}

\subsection{Full Formulation}

Given the new objective function, constraints, and decision variables, with the addition of orbital maneuvering as a task and the reconfiguration stages, the full \REOSSP is denoted as follows: \hypertarget{REOSSP}{}
\begin{equation}
    \begin{split}
        \max & \quad z_{R} \\
        \text{s.t.} & \quad \text{Constraints \eqref{REOSSP:dv}--\eqref{REOSSP:obj}, \eqref{REOSSP:flow_constraints}--\eqref{REOSSP:Battery_Constrain}}
    \end{split} \tag{\textsf{REOSSP}}
    \label{REOSSP:REOSSP}
\end{equation}

The objective function again balances the number of occurrences of data downlink to any ground stations $g \in \mathcal{G}$ with the less heavily weighted number of occurrences of the observation of any targets $p \in \mathcal{P}$. Simultaneously, the constraints restrict the problem to the observation of data prior to data downlink, the assurance of available data storage, the assurance of available battery power, and feasible constellation reconfigurability with consideration of a maximum fuel budget. Additionally, the \REOSSP provides the schedule of charging, observation, data downlink, and orbital maneuvers, as well as the path of orbital maneuvers, of a constellation of satellites over a given time horizon. Appendix~A contains a mathematical proof validating the identical nature of the \EOSSP and \REOSSP in the event that $c^k_{\max} = 0, \forall k \in \mathcal{K}$, wherein no satellite may perform orbital maneuvers.

\section{Rolling Horizon Procedure} \label{sec:RHP}

The RHP solution method is an iterative algorithmic approach to the REOSSP that divides the problem into sequential subproblems, each incorporating future lookahead stages and updating the problem state as the subproblems are solved at reconfiguration stages. This solution method is introduced to address the computational intractability that may occur in large instances of the \REOSSP and will be referred to as the \RHP for brevity. The \RHP analyzes the impact of decisions on future lookahead stages to make a more informed decision in the current stage with a much smaller problem scale. 

The \RHP heavily mirrors the \REOSSP in all aspects, with a few key differences through the expansion of time parameters to separate the problem into piecewise subproblems along the time horizon. As such, only the differences between the \RHP and \REOSSP will be explicitly stated. Each subproblem, denoted RHP($s,L$), consisted of a control stage $s \in \mathcal{S} \setminus \{0\}$ in which the solution is isolated and $L$ lookahead stages which aid in the decision-making process. Therefore, the \RHP reformats an $S$-stage \REOSSP to $S-L$ total subproblems, and the cumulative solution of the \RHP is composed of the solutions to each subproblem. Within each subproblem, the parameter $\ell \in \mathcal{L} = \{s,\ldots,s+L\}$ is used to denote the lookahead stage within a subproblem, replacing the $s$ used in the \REOSSP (e.g., ${V}^{\ell k}_{tjp}$ and ${c}^{\ell k}_{ij}$).

Certain parameters that are relevant throughout the entire schedule duration must be updated within each subproblem of the \RHP. Firstly, the maximum transfer cost budget of each subproblem depends on the orbital maneuvers performed within previous subproblems; therefore, the maximum transfer cost budget of each subproblem is given as $c^{sk}_{\max}$, which is obtained via Eq.~\eqref{RHP:cost_OE}, where $\tilde{x}^{nk}_{ij}$ is the optimal orbital maneuver decision variable from previous subproblems. Secondly, the initial satellite orbital slots of each control stage depend on the orbital maneuvers performed within previous subproblems; therefore, the set $\tilde{\mathcal{J}}^{s-1, k}$ is defined according to Eq.~\eqref{RHP:Jtilde_OE}. Finally, the initial data storage and battery storage values depend upon the tasks performed within previous subproblems, wherein these values will be assigned as $d^{sk}_1$ and $b^{sk}_1$ obtained via Eqs.~\eqref{RHP:d1_OE} and~\eqref{RHP:b1_OE}, respectively. In Eqs.~\eqref{RHP:d1_OE} and~\eqref{RHP:b1_OE}, the decision variables $\tilde{y}^{nk}_{tp}$, $\tilde{q}^{nk}_{tg}$, and $\tilde{h}^{nk}_t$ correspond to the optimal decision variables of target observation, ground station downlink, and solar charging within previous subproblems, respectively. With respect to the first stage $s=1$, all values are assumed to be the same as in the \REOSSP such that $c^{1k}_{\max} = c^k_{\max}, \forall k \in \mathcal{K}$, $d^{1k}_1 = D^k_{\min}, \forall k \in \mathcal{K}$, and $b^{1k}_1 = B^k_{\max}, \forall k \in \mathcal{K}$.
\begin{subequations}
    \begin{alignat}{2}
& c_{\max}^{sk} = c^k_{\max} - \sum^{s-1}_{n=1} \sum_{i \in \mathcal{J}^{n-1, k}} \sum_{j \in \mathcal{J}^{nk}} c^{nk}_{ij} \tilde{x}^{nk}_{ij}
\label{RHP:cost_OE} \\
& \tilde{\mathcal{J}}^{s-1, k} \coloneq \{j: \tilde{x}^{s-1, k}_{ij} = 1, i \in \mathcal{J}^{s-2, k}, j \in \mathcal{J}^{s-1, k} \} 
\label{RHP:Jtilde_OE} \\
& d^{sk}_1 = D^k_{\min} + \sum^{s-1}_{n=1} \sum_{t \in \mathcal{T}^n}\sum_{p \in \mathcal{P}} D_{\text{obs}} \tilde{y}^{nk}_{tp} - \sum^{s-1}_{n=1} \sum_{t \in \mathcal{T}^n} \sum_{g \in \mathcal{G}} D_{\text{comm}} \tilde{q}^{nk}_{tg}
\label{RHP:d1_OE} \\
\begin{split}
    & b^{sk}_{1} = B^k_{\text{max}} + \sum^{s-1}_{n=1} \sum_{t \in \mathcal{T}^n} B_{\text{charge}}\tilde{h}^{nk}_t - \sum^{s-1}_{n=1} \sum_{t \in \mathcal{T}^n} \sum_{p \in \mathcal{P}} B_{\text{obs}}\tilde{y}^{nk}_{tp} \\
    & \qquad - \sum^{s-1}_{n=1} \sum_{t \in \mathcal{T}^n} \sum_{g \in \mathcal{G}} B_{\text{comm}} \tilde{q}^{nk}_{tg} - \sum^{s-1}_{n=1} \sum_{i \in \mathcal{J}^{n-1,k}} \sum_{j \in \mathcal{J}^{nk}}B_{\text{recon}} \tilde{x}^{nk}_{ij} - (s-1)T^s B_{\text{time}}
\end{split} 
\label{RHP:b1_OE}
    \end{alignat}
    \label{RHP:OE}
\end{subequations}

\subsection{Decision Variables and Indicator Variables}

Each subproblem RHP($s,L$) contains the same decision and indicator variables as the \REOSSP, with additional consideration that $\mathcal{J}^{\ell - 1, k}$ is only a singleton set when $s = 1$.
\begin{subequations}
    \begin{alignat}{2}
        x^{\ell k}_{ij} & \in \{0, 1\}, \quad && \forall \ell \in \mathcal{L}, \forall k \in \mathcal{K}, \forall i \in \mathcal{J}^{\ell-1, k}, \forall j \in \mathcal{J}^{\ell k} \label{RHP:x} \\
        y^{\ell k}_{tp} & \in \{0, 1\}, \quad && \forall \ell \in \mathcal{L}, \forall t \in \mathcal{T}^\ell, \forall p \in \mathcal{P}, \forall k \in \mathcal{K} \label{RHP:y} \\
        q^{\ell k}_{tg} & \in \{0, 1\}, \quad && \forall \ell \in \mathcal{L}, \forall t \in \mathcal{T}^\ell, \forall g \in \mathcal{G}, \forall k \in \mathcal{K} \label{RHP:q} \\
        h^{\ell k}_t & \in \{0, 1\}, \quad && \forall \ell \in \mathcal{L}, \forall t \in \mathcal{T}^\ell, \forall k \in \mathcal{K} \label{RHP:h} \\
        d^{\ell k}_t & \in [D^k_{\min}, D^k_{\max}], \quad && \forall \ell \in \mathcal{L}, \forall t \in \mathcal{T}^\ell, \forall k \in \mathcal{K} \label{RHP:d} \\
        b^{\ell k}_t & \in [B^k_{\min}, B^k_{\max}], \quad && \forall \ell \in \mathcal{L}, \forall t \in \mathcal{T}^\ell, \forall k \in \mathcal{K} \label{RHP:b}
    \end{alignat}
    \label{RHP:dv}
\end{subequations}

\subsection{Objective Function}

As with the \REOSSP, the objective function of each subproblem RHP($s,L$) maximizes the number of observations and subsequent downlink occurrences of data to ground stations with weighting $C$ to prioritize the downlink of data. The objective function value of each subproblem is denoted as $z^s_{\text{RHP}}$ obtained via
\begin{equation}
    z^s_{\text{RHP}} = \sum_{k \in \mathcal{K}} \sum_{\ell \in \mathcal{L}} \sum_{t \in \mathcal{T}^{\ell}} \left( \sum_{g \in \mathcal{G}} C q^{\ell k}_{tg} + \sum_{p \in \mathcal{P}} y^{\ell k}_{tp} \right) 
    \label{RHP:obj}
\end{equation}
where $z^s_{\text{RHP}}$ is the objective function value of subproblem $s$. Additionally, a figure of merit $Z_{\text{RHP}}$ is defined for the \RHP at the end of this section relative to the optimal solution to each subproblem. The subscript $\text{RHP}$ is used to distinguish between the objective function values and figures of merit for the \EOSSP and \REOSSP. 

\subsection{Constraints}

Each subproblem RHP($s,L$) expands the functionality of the constraints of the \REOSSP according to the updating equations Eq.~\eqref{RHP:OE}, which update key parameters between subproblems.

\subsubsection{Orbital Maneuver Path Continuity Constraints}

The first set of constraints mirror constraints~\eqref{REOSSP:flow_constraints}, while constraints~\eqref{RHP:flow_from_initial} and~\eqref{RHP:flow_cost} depend on the updated parameters $\tilde{\mathcal{J}}^{s-1, k}$ and $c^{sk}_{\max}$, respectively.
\begin{subequations}
    \begin{alignat}{2}
& \sum_{j \in \mathcal{J}^{sk}} x^{sk}_{ij} = 1, \quad && \forall k \in \mathcal{K}, \forall i \in \tilde{\mathcal{J}}^{s-1, k} 
\label{RHP:flow_from_initial}\\ 
& \sum_{j \in \mathcal{J}^{\ell+1, k}} x^{\ell+1, k}_{ij} - \sum_{j' \in \mathcal{J}^{\ell-1, k}} x^{\ell k}_{j'i} = 0, \quad && \forall \ell \in \mathcal{L} \setminus \{s+L\}, \forall k \in \mathcal{K}, \forall i \in \mathcal{J}^{\ell k} \label{RHP:flow_from_0<ell<s+L-1}\\ 
& \sum_{\ell \in \mathcal{L}} \sum_{i \in \mathcal{J}^{\ell-1, k}} \sum_{j \in \mathcal{J}^{\ell k}} c^{\ell k}_{ij} x^{\ell k}_{ij} \leq c^{sk}_{\max}, \quad && \forall k \in \mathcal{K} 
\label{RHP:flow_cost}
    \end{alignat}
    \label{RHP:flow_constraints}
\end{subequations}

\subsubsection{Time Window Constraints}

The second set of constraints directly mirror constraints~\eqref{REOSSP:Visibility} relating to VTWs and task performance opportunities.
\begin{subequations}
    \begin{alignat}{2}
& \sum_{i \in \mathcal{J}^{\ell-1, k}} \sum_{j \in \mathcal{J}^{\ell k}} V^{\ell k}_{tjp} x^{\ell k}_{ij} \ge y^{\ell k}_{tp}, \quad && \forall \ell \in \mathcal{L}, \forall t \in \mathcal{T}^\ell, \forall p \in \mathcal{P}, \forall k \in \mathcal{K} 
\label{RHP:target_visibility}\\ 
& \sum_{i \in \mathcal{J}^{\ell-1, k}} \sum_{j \in \mathcal{J}^{\ell k}} W^{\ell k}_{tjg} x^{\ell k}_{ij} \ge q^{\ell k}_{tg}, \quad && \forall \ell \in \mathcal{L}, \forall t \in \mathcal{T}^\ell, \forall g \in \mathcal{G}, \forall k \in \mathcal{K} 
\label{RHP:gs_visibility}\\ 
& \sum_{i \in \mathcal{J}^{\ell-1, k}} \sum_{j \in \mathcal{J}^{\ell k}} H^{\ell k}_{tj} x^{\ell k}_{ij} \ge h^{\ell k}_t, \quad && \forall \ell \in \mathcal{L}, \forall t \in \mathcal{T}^\ell, \forall k \in \mathcal{K} 
\label{RHP:sun_visibility}\\  
& \sum_{p \in \mathcal{P}} y^{\ell k}_{tp} + \sum_{g \in \mathcal{G}} q^{\ell k}_{tg} + h^{\ell k}_t \leq 1, \quad && \forall \ell \in \mathcal{L}, \forall t \in \mathcal{T}^\ell, \forall k \in \mathcal{K} 
\label{RHP:obvs-down_overlap}
    \end{alignat}
    \label{RHP:visibility}
\end{subequations}

\subsubsection{Data Tracking and Storage Constraints}

The data tracking and storage constraints in each subproblem RHP($s,L$) closely mirror constraints~\eqref{REOSSP:Data}, with further breakdown pertaining to the consideration of the updated parameter $d^{sk}_1$ from Eq.~\eqref{RHP:d1_OE}. As such, constraints~\eqref{RHP:d-track_T1} account for the updated parameter, constraints~\eqref{RHP:d-track_s_stage} account for the rest of the $s$ stage, constraints~\eqref{RHP:d-track_l_stage} account for times within each other stage, and constraints~\eqref{RHP:d-track_next_stage} account for the stage gaps. Constraints~\eqref{RHP:d<max} and~\eqref{RHP:d>0} are unchanged apart from the consideration of $\ell$.

\begin{subequations}
    \begin{alignat}{2}
& d^{sk}_2 = d^{sk}_1 + \sum_{p \in \mathcal{P}} D_{\text{obs}} y^{sk}_{1p} - \sum_{g \in \mathcal{G}} D_{\text{comm}} q^{sk}_{1g}, \quad && \forall k \in \mathcal{K} 
\label{RHP:d-track_T1}\\
& d_{t+1}^{sk} = d^{sk}_t + \sum_{p \in \mathcal{P}} D_{\text{obs}} y^{sk}_{tp} - \sum_{g \in \mathcal{G}} D_{\text{comm}} q^{sk}_{tg}, \quad && \forall k \in \mathcal{K}, \forall t \in \mathcal{T}^s \setminus \{1,T^s\} 
\label{RHP:d-track_s_stage}\\ 
& d^{\ell k}_{t+1} = d^{\ell k}_t + \sum_{p \in \mathcal{P}} D_{\text{obs}} y^{\ell k}_{tp} - \sum_{g \in \mathcal{G}} D_{\text{comm}} q^{\ell k}_{tg}, \quad && \forall \ell \in \mathcal{L} \setminus \{s\}, \forall t \in \mathcal{T}^\ell \setminus \{T^\ell\}, \forall k \in \mathcal{K} 
\label{RHP:d-track_l_stage}\\
& d^{\ell+1,k}_1 = d^{\ell k}_{T^\ell} + \sum_{p \in \mathcal{P}} D_{\text{obs}} y^{\ell k}_{T^{\ell} p} - \sum_{g \in \mathcal{G}} D_{\text{comm}} q^{\ell k}_{T^{\ell} g}, \quad && \forall \ell \in \mathcal{L} \setminus \{s+L\}, \forall k \in \mathcal{K} 
\label{RHP:d-track_next_stage}\\ 
& d^{\ell k}_t + \sum_{p \in \mathcal{P}} D_{\text{obs}} y^{\ell k}_{tp} \leq D^k_{\max}, \quad && \forall \ell \in \mathcal{L} , \forall t \in \mathcal{T}^\ell , \forall k \in \mathcal{K} 
\label{RHP:d<max}\\ 
& d^{\ell k}_t - \sum_{g \in \mathcal{G}} D_{\text{comm}} q^{\ell k}_{tg} \ge D^k_{\min}, \quad && \forall \ell \in \mathcal{L}, \forall t \in \mathcal{T}^\ell , \forall k \in \mathcal{K} 
\label{RHP:d>0}
    \end{alignat}
    \label{RHP:Data}
\end{subequations}

\subsubsection{Battery Tracking Constraints}

The battery tracking constraints of each subproblem RHP($s,L$) also closely mirror constraints~\eqref{REOSSP:Battery_Track}, with similar breakdown relative to the updated parameter $b^{sk}_1$ from Eq.~\eqref{RHP:b1_OE}. Constraints~\eqref{RHP:b-track_T1} account for the updated parameter, constraints~\eqref{RHP:b-track_s_stage} account for the rest of the $s$ stage, constraints~\eqref{RHP:b-track_l_stage} account for times within all other stages, and constraints~\eqref{RHP:b-track_next_stage} account for the stage gaps.
\begin{subequations}
    \begin{alignat}{2}
& b^{sk}_2 = b^{sk}_1 + B_{\text{charge}} h^{sk}_1 - \sum_{p \in \mathcal{P}} B_{\text{obs}} y^{sk}_{1p} - \sum_{g \in \mathcal{G}} B_{\text{comm}} q^{sk}_{1g} - B_{\text{time}}, \quad && \forall k \in \mathcal{K} 
\label{RHP:b-track_T1}\\ 
& b^{sk}_{t+1} = b^{sk}_t +B_{\text{charge}} h^{sk}_t - \sum_{p \in \mathcal{P}} B_{\text{obs}} y^{sk}_{tp} - \sum_{g \in \mathcal{G}} B_{\text{comm}} q^{sk}_{tg} - B_{\text{time}}, \quad && \forall t \in \mathcal{T}^s \setminus \{1,T^s\}, \forall k \in \mathcal{K} 
\label{RHP:b-track_s_stage}\\ 
& b^{\ell k}_{t+1} = b^{\ell k}_t + B_{\text{charge}} h^{\ell k}_t - \sum_{p \in \mathcal{P}} B_{\text{obs}} y^{\ell k}_{tp} - \sum_{g \in \mathcal{G}} B_{\text{comm}} q^{\ell k}_{tg} - B_{\text{time}}, \quad && \forall \ell \in \mathcal{L} \setminus \{s\} , \forall t \in \mathcal{T}^\ell \setminus \{T^\ell\}, \forall k \in \mathcal{K} 
\label{RHP:b-track_l_stage} \\
\begin{split}
    & b^{\ell+1,k}_1 = b^{\ell k}_{T^\ell} + B_{\text{charge}} h^{\ell k}_{T^{\ell}} - \sum_{p \in \mathcal{P}} B_{\text{obs}} y^{\ell k}_{T^{\ell} p} - \sum_{g \in \mathcal{G}} B_{\text{comm}} q^{\ell k}_{T^{\ell} g} - \\
    & \qquad \sum_{i \in \mathcal{J}^{\ell k}} \sum_{j \in \mathcal{J}^{\ell+1, k}} B_{\text{recon}} x^{\ell+1, k}_{ij} - B_{\text{time}},
\end{split} \quad && \forall \ell \in \mathcal{L} \setminus \{s+L\} , \forall k \in \mathcal{K} 
\label{RHP:b-track_next_stage}
    \end{alignat}
    \label{RHP:Battery_Track}
\end{subequations}

\subsubsection{Battery Storage Constraints}

Finally, the battery storage constraints of each subproblem RHP($s,L$) directly mirror constraints~\eqref{REOSSP:Battery_Constrain} for the feasible restriction of battery capacities.
\begin{subequations}
    \begin{alignat}{2}
& b^{\ell k}_t + B_{\text{charge}} h^{\ell k}_t \leq B^k_{\max}, \quad && \forall \ell \in \mathcal{L}, \forall t \in \mathcal{T}^\ell, \forall k \in \mathcal{K} 
\label{RHP:b<max}\\ 
& b^{\ell k}_t - \sum_{p \in \mathcal{P}} B_{\text{obs}} y^{\ell k}_{tp} - \sum_{g \in \mathcal{G}} B_{\text{comm}} q^{\ell k}_{tg} - B_{\text{time}} \ge B^k_{\min}, \quad && \forall \ell \in \mathcal{L}, \forall t \in \mathcal{T}^\ell, \forall k \in \mathcal{K} 
\label{RHP:b>0}\\ 
\begin{split}
& b^{\ell k}_{T^\ell} - \sum_{p \in \mathcal{P}} B_{\text{obs}} y^{\ell k}_{T^{\ell} p} - \sum_{g \in \mathcal{G}} B_{\text{comm}} q^{\ell k}_{T^{\ell} g} - \\
& \qquad \sum_{i \in \mathcal{J}^{\ell k}}\sum_{j \in \mathcal{J}^{\ell+1, k}} B_{\text{recon}}x^{\ell+1,k}_{ij} - B_{\text{time}} \ge B^k_{\min},
\end{split} \quad && \forall \ell \in \mathcal{L} \setminus \{s+L\}, \forall k \in \mathcal{K} 
\label{RHP:b>0_end}\\ 
& b^{sk}_{1} - \sum_{i \in \mathcal{J}^{s-1,k}} \sum_{j \in \mathcal{J}^{s k}} B_{\text{recon}} x^{sk}_{ij} \ge B^k_{\min}, \quad && \forall k \in \mathcal{K} 
\label{RHP:b>0_gap}
    \end{alignat}
    \label{RHP:Battery_constraint}
\end{subequations}

\subsection{Full Formulation}

With the definition of the objective function, constraints, and decision variables, subproblem RHP($s,L$) is given as
\begin{equation}
    \begin{split}
        \max & \quad z^s_{\text{RHP}} \\
        \text{s.t.} & \quad \text{Constraints \eqref{RHP:dv}--\eqref{RHP:Battery_constraint}}
    \end{split} \tag{RHP($s,L$)}
    \label{RHP:optimization}
\end{equation}
Within each subproblem, the initial conditions of the control stage $s$ are updated using Eqs.~\eqref{RHP:OE}. The subproblem RHP($s,L$) is then solved iteratively for all control stages $s \in \{1, 2, \ldots, S-L\}$ through the use of the full \RHP algorithm, defined in Algorithm~\ref{RHP:algorithm}, where the objective of each subproblem RHP($s,L$) is to maximize the number of downlink occurrences to any ground station and the number of target observations, with a heavier weighting on the number of downlink occurrences. The constraints restrict the data and battery storage values of each satellite, as well as limit the propellant budget for orbital maneuvers. Upon obtainment of the solution of each subproblem, the optimal decision variables for the first stage of each subproblem are appended linearly within the time horizon to obtain the set of optimal decision variables throughout the entire time horizon, given as $(\tilde{x}^{sk}_{ij},\tilde{y}^{sk}_{tp},\tilde{q}^{sk}_{tg},\tilde{h}^{sk}_{t},\tilde{d}^{sk}_{t},\tilde{b}^{sk}_{t})$ in Algorithm~\ref{RHP:algorithm}. The \RHP provides the schedule of each satellite, including the charging, observing, downlinking, and orbital maneuvering sequence to provide a complete overview of the tasks performed by each satellite. The \RHP includes a figure of merit found using the optimal solution, corresponding to the amount of data downlinked to ground stations by the optimal decision variable $\tilde{q}^{sk}_{tg}$ and obtained via 

\begin{equation}
    Z_{\text{RHP}} = \sum_{k \in \mathcal{K}} \sum_{s \in \mathcal{S} \setminus \{0\}} \sum_{t \in \mathcal{T}^{s}} \sum_{g \in \mathcal{G}} D_{\text{comm}} \tilde{q}^{sk}_{tg} 
    \label{RHP:merit}
\end{equation}

\begin{algorithm}[!ht] 
    \hypertarget{RHP}{}
    \DontPrintSemicolon
    \caption{\textcolor{myblue}{\textsf{REOSSP-RHP}}}
    \label{RHP:algorithm}
    \For{$s = 1, \ldots, S-L-1$}{
        \If{$s \ge 2$}{
            Update $c^{sk}_{\max}, \tilde{\mathcal{J}}^{s-1, k}, d^{sk}_1,$ and $b^{sk}_1$ Using Eq.~\eqref{RHP:OE} \;
            } 
        Solve RHP($s,L$) and store: $z^s_{\text{RHP}}$ and 
        $(\tilde{x}^{sk}_{ij},\tilde{y}^{sk}_{tp},\tilde{q}^{sk}_{tg},\tilde{h}^{sk}_{t},\tilde{d}^{sk}_{t},\tilde{b}^{sk}_{t})$ \;
        } 
    $s \gets S-L$ \;
    Solve RHP($S-L,L$) and store: $\{z^1_{\text{RHP}},\ldots,z^S_{\text{RHP}}\}$ and 
    $(\tilde{x}^{sk}_{ij},\tilde{y}^{sk}_{tp},\tilde{q}^{sk}_{tg},\tilde{h}^{sk}_{t},\tilde{d}^{sk}_{t},\tilde{b}^{sk}_{t})\coloneq\{(\tilde{x}^{1k}_{ij},\tilde{y}^{1k}_{tp},\tilde{q}^{1k}_{tg},\tilde{h}^{1k}_{t},\tilde{d}^{1k}_{t},\tilde{b}^{1k}_{t}),\ldots,(\tilde{x}^{Sk}_{ij},\tilde{y}^{Sk}_{tp},\tilde{q}^{Sk}_{tg},\tilde{h}^{Sk}_{t},\tilde{d}^{Sk}_{t},\tilde{b}^{Sk}_{t})\}$ \;
    \Return $(\tilde{x}^{sk}_{ij},\tilde{y}^{sk}_{tp},\tilde{q}^{sk}_{tg},\tilde{h}^{sk}_{t},\tilde{d}^{sk}_{t},\tilde{b}^{sk}_{t})$\;
\end{algorithm}

\section{Computational Experiments} \label{sec:Experiments}

The \REOSSP, solved using a commercial software (therefore theoretically returning an exact optimal solution if given adequate computational resources and runtime) and referred to as the \REOSSPExact, and \RHP are benchmarked against the \EOSSP and compared to one another through a wide variety of schedule characteristics to evaluate the performance of each formulation relative to a baseline nonreconfigurable scheduling solution. Benchmarking experiments are conducted in two ways: including random instances of randomly varied schedule characteristics and a case study involving Hurricane Sandy, a highly dynamic natural disaster. The objective function values, $z_{E}$, $z_{R}$, and $z_{\text{RHP}}$ for the \EOSSP, \REOSSPExact, and \RHP, respectively, compare the performance of each scheduling problem. Additionally, the figures of merit, VTW of targets and ground stations, amount of data observed but not downlinked, amount of battery required to perform operations, and the orbital maneuver cost consumed by the satellites in the resulting schedules of the \REOSSPExact and \RHP are presented for analysis. Finally, the resultant schedule is used to gain insight into the order of task performance, especially concerning the case study. All computational experiments assume deterministic problems, where characteristics of a schedule are known \textit{a priori}.

All random instances and case studies are computed on a platform equipped with an Intel Core i9-12900 $2.4$ GHz (base frequency) CPU processor ($16$ cores and $24$ logical processors) and $64$ GB of RAM. While these specifications are not indicative of the capabilities of a satellite for onboard scheduling, they are relevant to demonstrate the computational efficiency of the formulations tested. Additionally, the \EOSSP, \REOSSPExact, and \RHP are programmed in MATLAB \cite{MATLAB} with the use of YALMIP \cite{YALMIP} and are solved using the commercial software package Gurobi Optimizer (version 12.0.0) with default settings apart from an assigned runtime.

\subsection{Random Instances} \label{subsec:random_experiments}

A set of $24$ random instances with random schedule characteristics demonstrates the capability of each scheduling problem in response to varied parameters. Each instance draws combinations of given parameters, creating variation in the number of stages, satellites, and orbital slot options available. The sets for such combinations are $S \in \{8, 9, 12\}$, $K \in \{5, 6\}$, and $J^{sk} \in \{20, 40, 60, 80\}, \forall s \in \mathcal{S} \setminus \{0\}, \forall k \in \mathcal{K}$, respectively. The Gurobi runtime limit for the \EOSSP, \REOSSPExact, and each subproblem of the \RHP is set to $60$ min (\SI{3600}{s}).

\subsubsection{Design of Random Instances}

Each of the $24$ random instances is unique to ensure that no instances have identical parameters, providing variation in ground station position, satellite orbits, and target position such that each scheduling problem is compared under a wide spectrum of conditions, thus not limiting the evaluation to a possibly biased set of conditions and providing a dynamic target through varied available rewards. All $K$ satellites are assigned inclined circular orbits with random orbital elements within specified ranges, such as an altitude between $600$ and \SI{1200}{km}, inclination between $40$ and \SI{80}{deg}, right ascension of ascending node (RAAN) between $0$ and \SI{360}{deg}, and argument of latitude between $0$ and \SI{360}{deg}. Additionally, all $J^{sk}$ orbital slots vary only in the argument of latitude (herein referred to as \textit{phase}) and are equally spaced between $0$ and \SI{360}{deg} with the inclusion of the initial phase. The number of ground stations is set as two distinct ground stations ($G=2$) and each is randomly assigned a latitude between \SI{80}{deg} South and \SI{80}{deg} North and a longitude between \SI{180}{deg} West and \SI{180}{deg} East, with an extra condition specifying that the ground station must be on land. The number of targets is set as $12$ distinct targets ($P=12$) where each target position varies on the same latitude and longitude range as the ground stations, with allowance of targets on water. Each target is assigned a binary masking within the visibility matrix $V^k_{tp}$ and $V^{sk}_{tjp}$ such that each target only has a visibility value of one on a time interval of $T/P$, thus allowing exactly one target to be considered as available at each time step. 

Certain parameters are fixed in each of the $24$ random instances, where all fixed parameters are included in Table~\ref{tab:Fixed_Params}, with necessary explanation following. Both $D_{\text{obs}}$ and $D_{\text{comm}}$ are converted from values in Mbps to discrete values in MB corresponding to a single time step. Additionally, $D_{\text{obs}}$ is the average of the two observation rates from Ref.~\cite{MODIS_specs}, and $D_{\text{comm}}$ is the X-band frequency of the SSTL-300 S1 spacecraft capable of Binary Phase-shift keying modulation from Ref.~\cite{DMC_specs}. Due to the small difference between $D_{\text{obs}}$ and $D_{\text{comm}}$, not all data observed in a single occurrence can be downlinked in a single occurrence, thus leaving a small amount of data in storage at the end of the schedule horizon. Furthermore, $B_{\text{obs}}$, $B_{\text{comm}}$, $B_{\text{recon}}$, $B_{\text{charge}}$, and $B_{\text{time}}$ are converted from kW to kJ, similarly to $D_{\text{obs}}$ and $D_{\text{comm}}$. Additionally, the value of $B^k_{\max}$ is converted from a $15$ A-h Battery using the average between $28$ and $33$ V as obtained from Ref.~\cite{DMC_specs}, and the value of $B_{\text{charge}}$ is generated using the specification of $\SI{2.44}{m^2}$ gallium arsenide solar cells from Ref.~\cite{DMC_specs} for power generation at $160$ to $180$ W/$\SI{}{m^2}$ \cite{Bradhorst_solarcells} in which the average of $170$ W is used.

\begin{table}[!ht]
    \centering
    \caption{Fixed random instance parameters}
    \begin{tabular}{ l c c  c c c }
        \hline \hline 
        Parameter & Value & Units & Parameter & Value & Units \\
        \hline 
        Schedule start time                     & January 1st, 2025, 00:00:00, UTC & --- & $T_r$                                   & $14$                     & days           \\
        $\Delta t$                              & $100$                            & s   & $T$                                     & \num{12096}              & ---            \\
        $D^k_{\min}, \forall k \in \mathcal{K}$ & $0$                              & GB  & $D^k_{\max}, \forall k \in \mathcal{K}$ & $128$ \cite{DMC_specs}   & GB             \\
        $D_{\text{obs}}$                        & $102.50$ \cite{MODIS_specs}      & MB  & $D_{\text{comm}}$                       & $100$ \cite{DMC_specs}   & MB             \\
        $B^k_{\min}, \forall k \in \mathcal{K}$ & $0$                              & kJ  & $B^k_{\max}, \forall k \in \mathcal{K}$ & $1647$ \cite{DMC_specs}  & kJ             \\
        $B_{\text{obs}}$                        & $16.26$ \cite{MODIS_specs}       & kJ  & $B_{\text{comm}}$                       & $1.20$ \cite{DMC_specs}  & kJ             \\
        $B_{\text{recon}}$                      & $0.50$ \cite{DMC_specs}          & kJ  & $B_{\text{charge}}$                     & $41.48$ \cite{DMC_specs} & kJ             \\
        $B_{\text{time}}$                       & $2$ \cite{DMC_specs}             & kJ  & $c^k_{\max}, \forall k \in \mathcal{K}$ & $750$                    & m/s            \\
        $C$                                     & $2$                              & --- & $L$                                     & $1$                      & ---            \\
        \hline \hline
    \end{tabular}
    \label{tab:Fixed_Params}
\end{table} 

In addition to the parameters listed, key methods are used to generate $V$, $W$, $H$, and $c^{sk}_{ij}$. Visibility matrices such as $V$ and $G$ are generated through the use of the Aerospace Toolbox \cite{MATLAB} by propagating satellite orbits and orbital slots, assigning each satellite a conical field of view (FOV) of \SI{45}{deg} for target observation, a conical communication range of \SI{120}{deg} for ground station communication, and evaluating the \texttt{access} function. This \texttt{access} function returns a binary value representing whether a given target object (either a target or ground station) is within the satellite FOV at any time step over the schedule duration. Separately, the visibility matrix $H$ is generated through the same propagation and the use of the \texttt{eclipse} function. This \texttt{eclipse} function returns a value between zero and one representing sunlight exposure via the dual cone method, with direct sunlight assigned one, penumbra assigned between zero and one, and umbra assigned zero \cite{Woodburn_2000}. The propagator used is the Simplified General Perturbation 4 (SGP4) model. Finally, the transfer cost $c^{sk}_{ij}$ is computed through the use of the circular coplanar phasing problem found in Ref.~\cite{Vallado2022}, in which orbital slot $i \in \mathcal{J}^{s-1, k}$ and $j \in \mathcal{J}^{sk}$ are defined as the boundary conditions. 

It should be noted that the chosen parameters for the random instances are specially selected to demonstrate the use and results of the formulated scheduling problems. Parameters such as constellation formation, satellite FOV, communication range, maximum transfer budget, and others, can be chosen to model desired EO applications. 

\subsubsection{Random Instance Results}

Table~\ref{tab:Random_results} contains the results of all $24$ random instances with associated specified varied parameters and ID numbers. The results include the optimal objective function values (or sum of all subproblem objective function values for the \RHP), the figures of merit, and the computational runtime for the \EOSSP, \REOSSPExact, and \RHP. The percent improvement of one schedule resultant objective function value over another is designated by $\gamma$. Additionally, the results include the total propellant used by the \REOSSPExact and \RHP in km/s. Finally, the results include the minimum, maximum, mean, and standard deviation values of the runtime, the total propellant used, and $\gamma$ values, where the reported values of the \REOSSPExact omit the instances that did not result in a feasible solution. 

\begin{landscape}

\begin{table}[!ht]
    \centering
    \caption{Random instance results}
    \fontsize{7.25}{10} \selectfont
    \begin{threeparttable} 
        \begin{tabular}{ l l l l  r r r  r r r r r  r r r r r r }
\hline \hline 
\multicolumn{4}{c}{Instance} & \multicolumn{3}{c}{\EOSSP} & \multicolumn{5}{c}{\REOSSPExact} & \multicolumn{6}{c}{\RHP} \\
\cmidrule(lr){1-4} \cmidrule(lr){5-7} \cmidrule(lr){8-12} \cmidrule(lr){13-18}
ID & $S$ & $K$ & $J^{sk}$ & $z_{E}$ & $Z_{E}$, GB & Runtime, min & $z_{R}$ & $Z_{R}$, GB & Runtime, min & $\sum c^{sk}_{ij}$, km/s & $\gamma$,\% & $\sum_{s = 1}^{S} z^s_{\text{RHP}}$ & $Z_{\text{RHP}}$, GB & Runtime, min & $\sum c^{\ell k}_{ij}$, km/s & $\gamma$,\% & $\gamma$,\% \\
~ & ~ & ~ & ~ & ~ & ~ & ~ & ~ & ~ & ~ & Total used & Over   & ~ & ~ & ~ & Total used & Over            & Over   \\
~ & ~ & ~ & ~ & ~ & ~ & ~ & ~ & ~ & ~ & propellant & \EOSSP & ~ & ~ & ~ & propellant & \REOSSP         & \EOSSP \\ 
~ & ~ & ~ & ~ & ~ & ~ & ~ & ~ & ~ & ~ & ~          & ~      & ~ & ~ & ~ & ~          & \textsf{-Exact} & ~      \\ 
\hline 
1      & 8      & 5      & 20     & 102    & 3.40   & 1.01   & 183        & 6.10  & 4.55         & 3.51 & 79.41  & 162    & 5.40   & 1.91   & 3.68   & -11.48 & 58.82  \\
2      & 8      & 5      & 40     & 144    & 4.80   & 1.21   & 288        & 9.60  & 27.16        & 2.96 & 100.00 & 252    & 8.40   & 5.02   & 3.63   & -12.50 & 75.00  \\
3      & 8      & 5      & 60     & 213    & 7.10   & 3.79   &---\tnote{a}& ---   & 60.15        & ---  & ---    & 393    & 13.10  & 18.24  & 3.66   & ---    & 84.51  \\
4      & 8      & 5      & 80     & 246    & 8.20   & 5.86   & ---        & ---   &---\tnote{b}  & ---  & ---    & 393    & 13.10  & 54.45  & 3.71   & ---    & 59.76  \\
5      & 8      & 6      & 20     & 193    & 6.40   & 7.68   & ---        & ---   & 60.12        & ---  & ---    & 325    & 10.80  & 14.40  & 4.35   & ---    & 68.75  \\
6      & 8      & 6      & 40     & 192    & 6.40   & 9.05   & 381        & 12.70 & 33.53        & 3.77 & 98.44  & 325    & 10.80  & 15.14  & 4.44   & -14.70 & 68.75  \\
7      & 8      & 6      & 60     & 246    & 8.20   & 11.90  & ---        & ---   & 60.21        & ---  & ---    & 470    & 15.70  & 22.81  & 4.43   & ---    & 91.46  \\
8      & 8      & 6      & 80     & 198    & 6.60   & 12.57  & ---        & ---   & 60.65        & ---  & ---    & 333    & 11.10  & 67.61  & 4.46   & ---    & 68.18  \\
9      & 9      & 5      & 20     & 129    & 4.30   & 2.80   & 249        & 8.30  & 21.13        & 3.10 & 93.02  & 225    & 7.50   & 13.79  & 3.57   & -9.64  & 74.42  \\
10     & 9      & 5      & 40     & 176    & 5.80   & 8.48   & ---        & ---   & 60.19        & ---  & ---    & 328    & 10.90  & 31.06  & 3.70   & ---    & 87.93  \\
11     & 9      & 5      & 60     & 162    & 5.40   & 5.71   & 327        & 10.90 & 22.84        & 3.19 & 101.85 & 297    & 9.90   & 8.26   & 3.68   & -9.17  & 83.33  \\
12     & 9      & 5      & 80     & 159    & 5.30   & 7.17   & ---        & ---   & 60.82        & ---  & ---    & 312    & 10.40  & 49.38  & 3.71   & ---    & 96.23  \\
13     & 9      & 6      & 20     & 252    & 8.40   & 1.62   & 440        & 14.70 & 14.33        & 4.16 & 75.00  & 366    & 12.20  & 14.43  & 4.26   & -16.82 & 45.24  \\
14     & 9      & 6      & 40     & 234    & 7.80   & 4.50   & 498        & 16.60 & 40.95        & 3.51 & 112.82 & 432    & 14.40  & 21.00  & 4.38   & -13.25 & 84.62  \\
15     & 9      & 6      & 60     & 115    & 3.80   & 3.72   & ---        & ---   & 60.17        & ---  & ---    & 235    & 7.80   & 12.68  & 4.47   & ---    & 105.26 \\
16     & 9      & 6      & 80     & 189    & 6.30   & 7.16   & ---        & ---   &---\tnote{b}  & ---  & ---    & 330    & 11.00  & 58.32  & 4.37   & ---    & 74.60  \\
17     & 12     & 5      & 20     & 150    & 5.00   & 1.06   & 324        & 10.80 & 9.64         & 2.79 & 116.00 & 279    & 9.30   & 5.03   & 3.63   & -13.89 & 86.00  \\
18     & 12     & 5      & 40     & 174    & 5.80   & 2.74   & 411        & 13.70 &60.10\tnote{c}& 3.15 & 136.21 & 372    & 12.40  & 9.64   & 3.68   & -9.49  & 113.79 \\
19     & 12     & 5      & 60     & 141    & 4.70   & 3.35   & ---        & ---   & 60.20        & ---  & ---    & 231    & 7.70   & 8.88   & 3.73   & ---    & 63.83  \\
20     & 12     & 5      & 80     & 213    & 7.10   & 5.71   & ---        & ---   &---\tnote{b}  & ---  & ---    & 381    & 12.70  & 54.64  & 3.73   & ---    & 78.87  \\
21     & 12     & 6      & 20     & 171    & 5.70   & 4.37   & 351        & 11.70 &60.11\tnote{c}& 4.01 & 105.26 & 270    & 9.00   & 5.04   & 4.43   & -23.08 & 57.89  \\
22     & 12     & 6      & 40     & 162    & 5.40   & 6.18   & ---        & ---   & 60.24        & ---  & ---    & 352    & 11.70  & 5.33   & 4.43   & ---    & 116.67 \\
23     & 12     & 6      & 60     & 260    & 8.60   & 5.95   & ---        & ---   & 60.45        & ---  & ---    & 398    & 13.20  & 14.04  & 4.43   & ---    & 53.49  \\
24     & 12     & 6      & 80     & 264    & 8.80   & 9.07   & ---        & ---   &---\tnote{b}  & ---  & ---    & 464    & 15.50  & 34.62  & 4.44   & ---    & 76.14  \\
\multicolumn{4}{c}{Minimum}            & --- & ---  & 1.01   & --- & ---          & 4.55   & 2.79   & 75.00  & --- & ---           & 1.91   & 3.57   & -23.08 & 45.24  \\
\multicolumn{4}{c}{Maximum}            & --- & ---  & 12.57  & --- & ---          & 60.11  & 4.16   & 136.21 & --- & ---           & 67.61  & 4.47   & -9.17  & 116.67 \\
\multicolumn{4}{c}{Mean}               & --- & ---  & 5.53   & --- & ---          & 29.43  & 3.41   & 101.80 & --- & ---           & 22.74  & 4.04   & -13.40 & 78.06  \\
\multicolumn{4}{c}{Standard deviation} & --- & ---  & 3.21   & --- & ---          & 19.40  & 0.46   & 17.73  & --- & ---           & 19.73  & 0.38   & 4.20   & 18.17  \\
\hline \hline
        \end{tabular}
        \begin{tablenotes}
           \item[a] A dash (-) indicates the trigger of the runtime limit of \SI{3600}{s} without an obtained feasible solution.
           \item[b] Gurobi error code 10001 (out of memory) occurred, and no feasible solution was obtained.
           \item[c] Although the runtime limit of \SI{3600}{s} was reached, a feasible solution was obtained. IDs $18$ and $21$ include a mixed integer programming gap of \SI{3.65}{\%} and \SI{0.28}{\%}, respectively.
        \end{tablenotes}
    \end{threeparttable}
    \label{tab:Random_results}
\end{table}
    
\end{landscape}

\begin{figure}[!ht]
    \centering
    \begin{subfigure}[h]{0.49\textwidth}
        \centering
        \includegraphics[width = \textwidth]{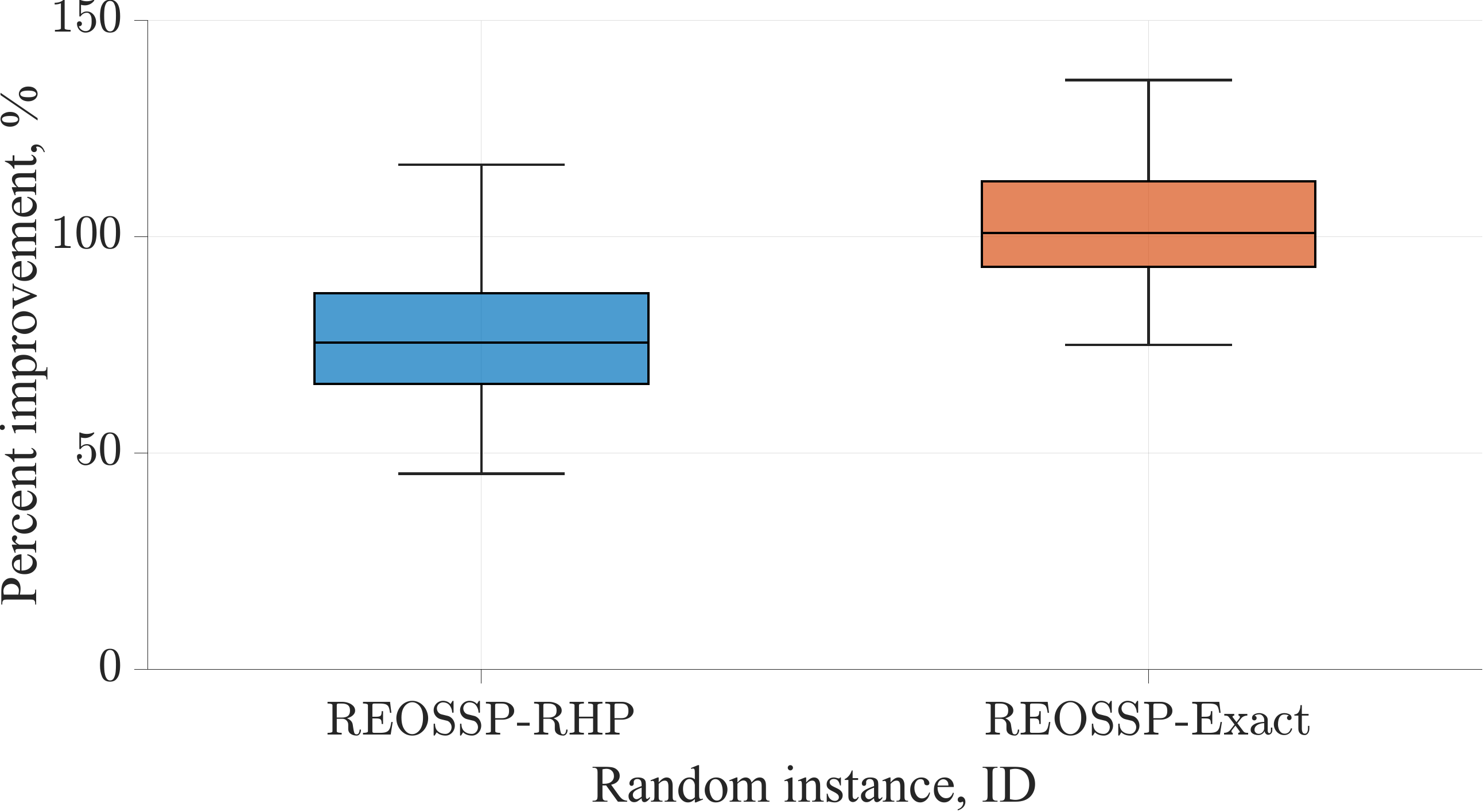}
        \caption{Percent of the \REOSSPExact and \RHP over the \EOSSP}
        \label{fig:random_percent_box}
    \end{subfigure}
    \hfill
    \begin{subfigure}[h]{0.49\textwidth}
        \centering
        \includegraphics[width = \textwidth]{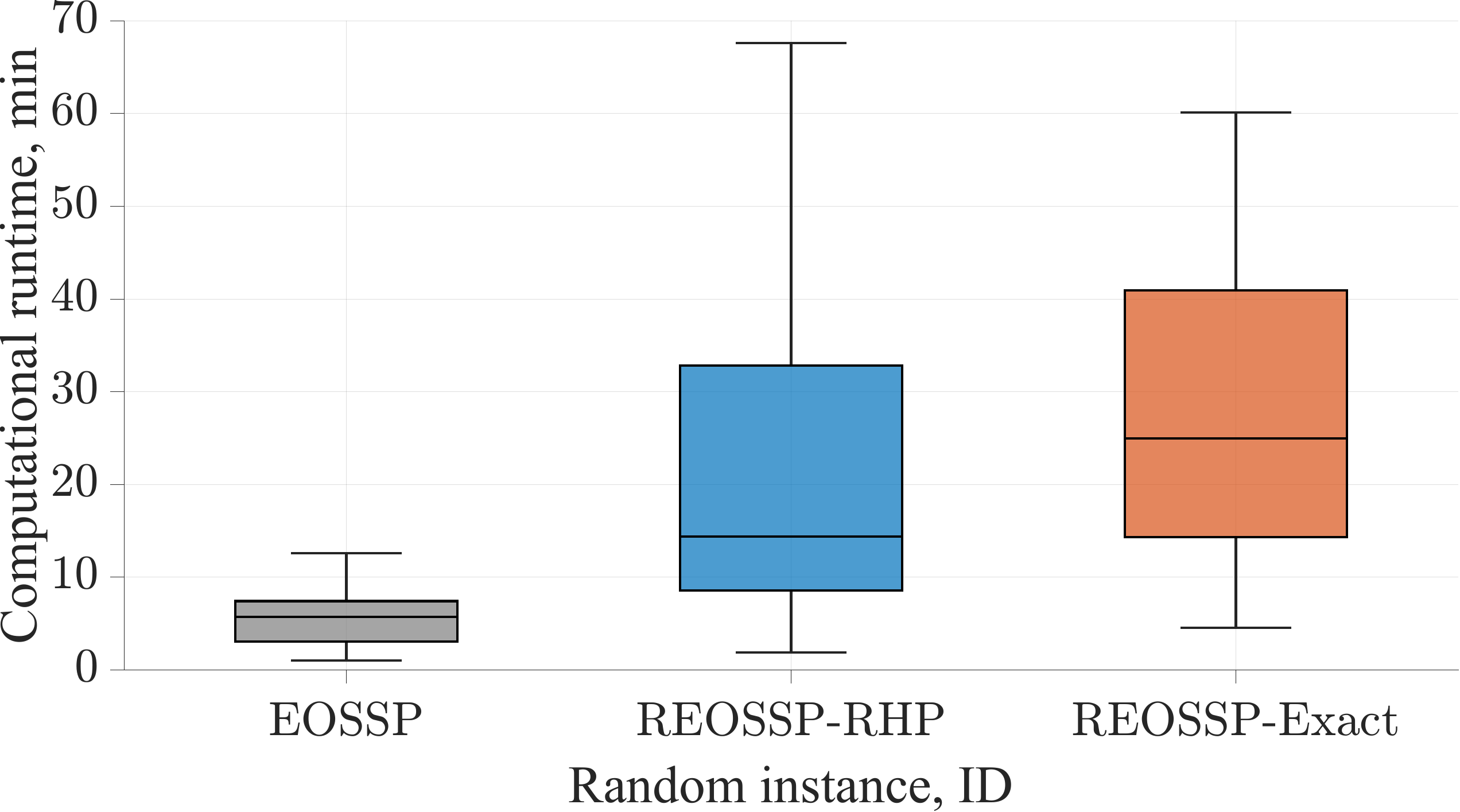}
        \caption{Runtime of the \EOSSP, \REOSSPExact, and \RHP}
        \label{fig:random_runtime_box}
    \end{subfigure}
    \caption{Results of the random instances.}
    \label{fig:rand_results_box}
\end{figure}

Both the \REOSSPExact and \RHP outperform the \EOSSP in every instance, apart from IDs $3, 5, 7, 8, 10, 12, 15, 19, 22, \text{and }23$ in which the \REOSSPExact triggers the runtime limit of \SI{3600}{s} without obtaining a feasible solution and IDs $4, 16, 20, \text{and }24$ in which the \REOSSPExact triggers Gurobi error code 10001 (out of memory), while the \RHP still outperforms the \EOSSP in these instances. Additionally, there are no instances in which the \RHP outperforms the \REOSSPExact, apart from those in which the \REOSSPExact does not obtain a feasible solution. This is expected since the \REOSSPExact has access to the entire schedule horizon when optimizing, while the \RHP only has access to a significantly reduced portion of the schedule. However, the \RHP performs computations significantly faster than the \REOSSPExact as a result of the significantly reduced problem scale within each subproblem of the \RHP, requiring less computation runtime in all instances apart from ID $13$, where the \RHP required approximately an additional \SI{6}{s} overall. The percent improvement of the \REOSSPExact and \RHP over the \EOSSP, as well as the runtime of the \EOSSP, \REOSSPExact, and \RHP, is shown as box charts in Fig.~\ref{fig:random_percent_box} and~\ref{fig:random_runtime_box}, respectively. Similar figures depicting the metrics of each individual instance are located in Appendix~B.

Statistically, the \REOSSPExact outperforms the \EOSSP with an average improvement of \SI{101.80}{\%} with a standard deviation of \SI{17.73}{\%}, a minimum improvement of \SI{0.00}{\%} (IDs $3$--$5, 7, 8, 10, 12, 15, 16, 19, 20, 22$--$24$), and a maximum improvement of \SI{136.21}{\%} (ID $18$). Similarly, the \RHP outperforms the \EOSSP with an average improvement of \SI{78.06}{\%} with a standard deviation of \SI{18.17}{\%}, a minimum improvement of \SI{45.24}{\%} (ID $13$), and a maximum improvement of \SI{116.67}{\%} (ID $22$). Furthermore, while the \RHP performs worse than the \REOSSPExact by an average of \SI{13.40}{\%} with a standard deviation of \SI{4.20}{\%}, the \RHP attains significant improvement over the \EOSSP in the event that the \REOSSPExact does not obtain a feasible solution within the time limit, during which the \RHP improves upon the \EOSSP by an average of \SI{80.41}{\%}. Additionally, the \RHP performs better relative to the \REOSSPExact when the number of stages is lower as a result of the larger $T^s$ value and thus larger subproblem size providing more problem information in each subproblem, while the \RHP performs worse as the number of stages increases. 

Additionally, the runtime of the \REOSSPExact is much greater than that of the \EOSSP, while the average runtime of the \RHP is slightly more than one-third of the \REOSSPExact runtime limit of \SI{3600}{s}. Contrarily, the runtime of the \RHP is greater than the runtime of the \REOSSPExact in ID $13$ and the runtime of the \RHP is greater than the runtime limit of the \REOSSPExact in ID $8$. Statistically, the runtime of the \REOSSPExact is \SI{29.43}{min} on average with a standard deviation of \SI{19.40}{min}, a minimum of \SI{4.55}{min} (ID $1$), and a maximum beyond the runtime limit in each instance with no feasible solution obtained as well as IDs $18$ and $21$, which obtained a suboptimal feasible solution after triggering the runtime limit. Separately, the runtime of the \RHP is lower at \SI{22.74}{min} on average with a standard deviation of \SI{19.73}{min}, a minimum of \SI{1.91}{min} (ID $1$), and a maximum of \SI{67.61}{min} (ID $8$). Additionally, the minimum runtime instance of the \EOSSP, \REOSSPExact, and \RHP is ID $1$, wherein the lowest combination of $K$ and $J$ are present, while the maximum runtime instances of each scheduling problem occur in the event that the larger values of $K$ and $J$ are present. As for the maximum runtime instance of the \RHP (ID $8$), the highest combination of $K$ and $J$ are present, while $S$ is the lowest possible value, thus resulting in the largest $T^s$ value and thus the largest subproblem size. As for those IDs where the \REOSSPExact does not obtain a feasible solution, the cause is twofold. The first cause originates in the complexity of the problem geometry, wherein satellite visibility overlap may be complicated. Separately, the second cause is due to the combinatorial nature of the problem; the solution space expands with the increase in decision variables, which is heavily influenced by the number of satellites and orbital slot options.

Finally, the main limitation of constellation reconfigurability is the fuel cost required to perform orbital maneuvers throughout the schedule horizon, wherein the optimal schedule of the \RHP uses more propellant than that of the \REOSSPExact. The average propellant used by the \REOSSPExact is \SI{3.41}{km/s} with a standard deviation of \SI{0.46}{km/s} (with the removal of the instances where no feasible solution was found). Meanwhile, the \RHP is more expensive, where the average propellant used is \SI{4.04}{km/s} with a slightly lower standard deviation of \SI{0.38}{km/s}. The \RHP is more expensive than the \REOSSPExact on average as a result of the limited information presented to each subproblem, potentially leading to slightly less optimal orbital maneuvers taking place.  

Overall, the average performance increase over the \EOSSP by both the \REOSSPExact and \RHP shows that each solution method can greatly outperform the \EOSSP despite the limited fuel budget for orbital maneuvers and additional complexity. Additionally, the value of the \RHP is demonstrated through the attainment of a feasible solution in a reduced amount of time on average, especially in those cases where the \REOSSPExact did not obtain a feasible solution. Furthermore, the \RHP outperforms the \REOSSPExact on average on account of the \REOSSPExact instances that did not result in a feasible solution. 

\subsection{Case Study - Hurricane Sandy} \label{subsec:case_study}

In addition to the random instances, Hurricane Sandy is selected for an in-depth analysis to provide real-world data to the scheduling problems. Hurricane Sandy occurred in 2012 and struck many islands throughout the Caribbean Sea as well as the Northeastern coast of the United States, causing $233$ deaths \cite{Sandy_information} and \$$65$ billion in damage \cite{Sandy_cost}. Hurricane Sandy achieved Category Three status on October $25$ at 6:00 AM with $115$ mph winds at the peak \cite{WUnderground}. 

\subsubsection{Design of Case Study}

\begin{figure}[!ht]
    \centering
    \begin{subfigure}[h]{0.5\textwidth}
        \centering
        \includegraphics[width=\textwidth]{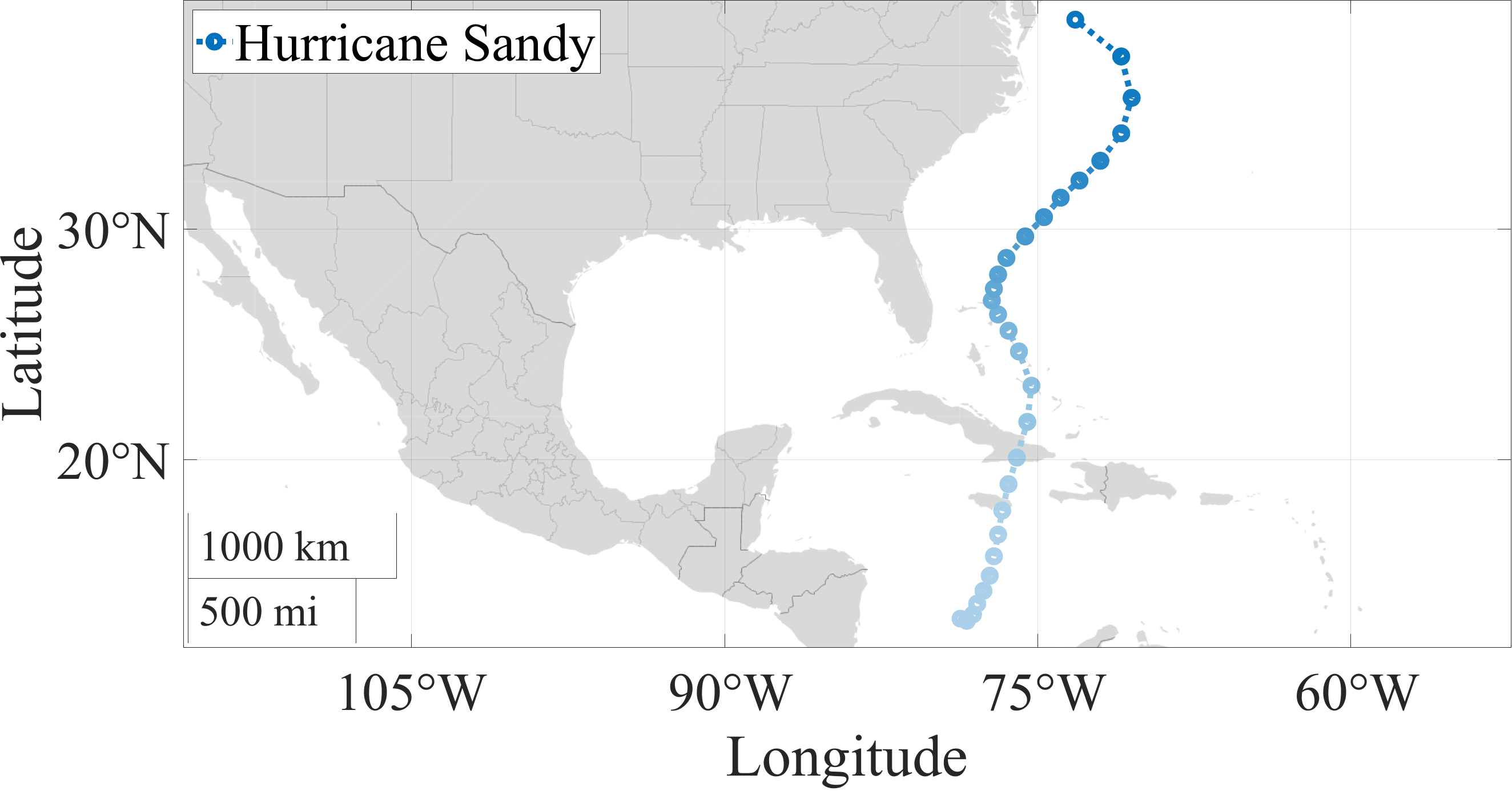}
        \caption{Path of Hurricane Sandy}
        \label{fig:Historical_tempests}
    \end{subfigure}
    \hfill
    \begin{subfigure}[h]{0.4\textwidth}
        \centering
        \includegraphics[width=0.825\textwidth]{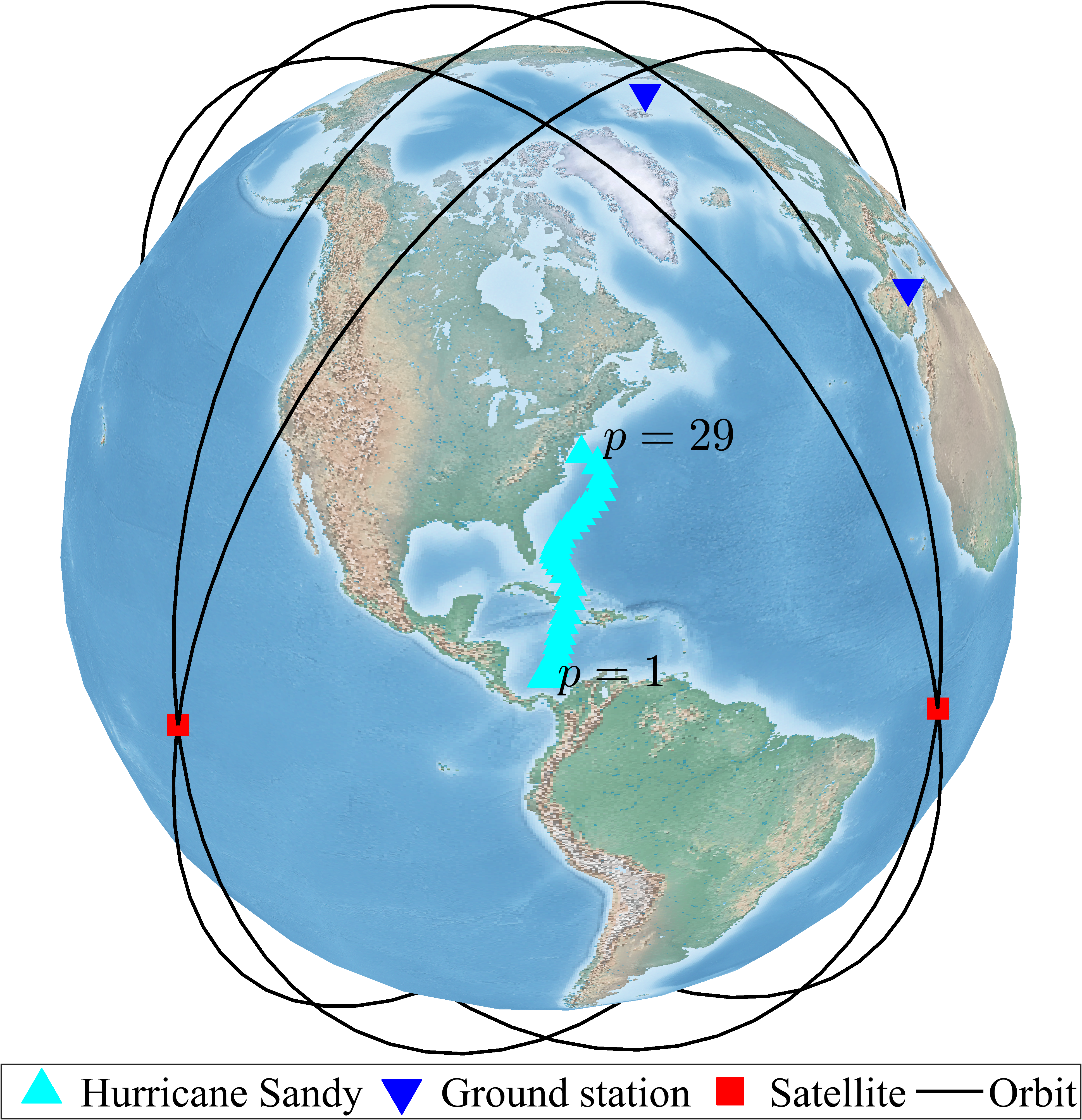}
        \caption{Three-dimensional locations of Hurricane Sandy, ground stations, and satellites}
        \label{fig:CaseStudyParameters}
    \end{subfigure}
    \caption{Locations of case study parameters.}
\end{figure}

For the case study, fixed schedule characteristics are selected rather than random characteristics to demonstrate realistic utilization of the scheduling problems further. Firstly, Hurricane Sandy is depicted in Fig.~\ref{fig:Historical_tempests} from the initial to the final occurrence of tropical storm status or higher (defined as containing winds between $39$ and $73$ mph \cite{TempestDefinitions_2022}). The historical path data utilized is gathered from Ref.~\cite{WUnderground} and tracks the eye of the storm at an interval of $6$ hours, and $T_r$ is set to \SI{7.25}{days}. Considering the tracking conditions, the number of targets (points along the storm path) for Hurricane Sandy is $P=29$. Additionally, binary masking is applied to $V^k_{tp}$ for each target in chronological order, such that target $p$ has a masking of one at time steps $t \in [1 + (p-1)(T/P), p(T/P)]$ and a masking of zero otherwise. This masking condition similarly applies for $V^{sk}_{tjp}$ and $V^{\ell k}_{tjp}$ along each stage in the time horizon. The same time step size, $\Delta t=\SI{100}{s}$, is used in the case study, resulting in a total number of time steps of $T=6264$. Additionally, ground stations utilized by disaster monitoring satellites and satellite orbits suitable for EO are employed. The ground stations are assigned as two used by the Disaster Monitoring Constellation, the first being the Svalbard Satellite Station located in Sweden at \SI{78.23}{deg} North and \SI{15.41}{deg} East and the second being a satellite station located in Boecillo Spain at \SI{41.54}{deg} North and \SI{4.70}{deg} West \cite{DMC_gs}. Furthermore, a Walker-delta constellation is used as it has become a common choice in recent years for satellite configurations. The Walker-delta constellation for $K=4$ satellites is defined at an altitude of \SI{709}{km} and takes the form $\SI{98.18}{deg}:4/4/0$, where the inclination is $\SI{98.18}{deg}$, and four satellites are located in four orbital planes with zero relative phasing between satellites. The constellation formation is selected to remain consistent with an identical FOV case in Ref.~\cite{Nag2015}. The three-dimensional locations of Hurricane Sandy over time, the two selected ground stations, and the initial satellite orbits are shown in Fig.~\ref{fig:CaseStudyParameters}. 

Secondly, additional orbital slot options that extend to changes in the orbital plane and phase are provided. Such orbital slots allow equally spaced changes in inclination, RAAN, and phase, in which the plane-change slots extend in the positive and negative direction of the initial condition. The resultant number of orbital slot options is $J^{sk} = f(2m - 1)$, where $f$ is the number of phasing options and $m$ is the number of plane-change options in each plane of either inclination or RAAN. The case studies make use of $f = 15$ and $m = 5$, resulting in $J^{sk} = 135, \forall s \in \mathcal{S} \setminus \{0\}, \forall k \in \mathcal{K}$. A demonstration of the orbital slot option space is shown in Fig.~8 of Ref.~\cite{Pearl2025Developing}. The maximum degree of separation between the initial condition and the furthest plane-change option is determined through the use of rearranged boundary value problems (assuming the entire budget, $c^k_{\max}$, is consumed) from Ref.~\cite{Vallado2022} with an applied scaling factor of \SI{75}{\%} to ensure the feasibility of multiple maneuvers. 

As a result of the increase in orbital slot options, additional cost computation algorithms are employed for all possible combinations of transfer types. This includes 1) phasing only, 2) inclination change only, 3) RAAN change only, 4) simultaneous inclination and RAAN change, 5) inclination change followed by phasing, 6) RAAN change followed by phasing, and 7) simultaneous inclination and RAAN change followed by phasing. As with phasing-only computations, analytical algorithms from Ref.~\cite{Vallado2022} are used. Finally, $S=8$ stages are used in the \REOSSPExact and the \RHP. All other parameters and parameter generation methods are kept identical to those in Sec.~\ref{subsec:random_experiments}. 

\subsubsection{Case Study Results}

\begin{figure}[!ht]
    \centering
    \begin{subfigure}[h]{0.49\textwidth}
        \centering
        \includegraphics[width = \textwidth]{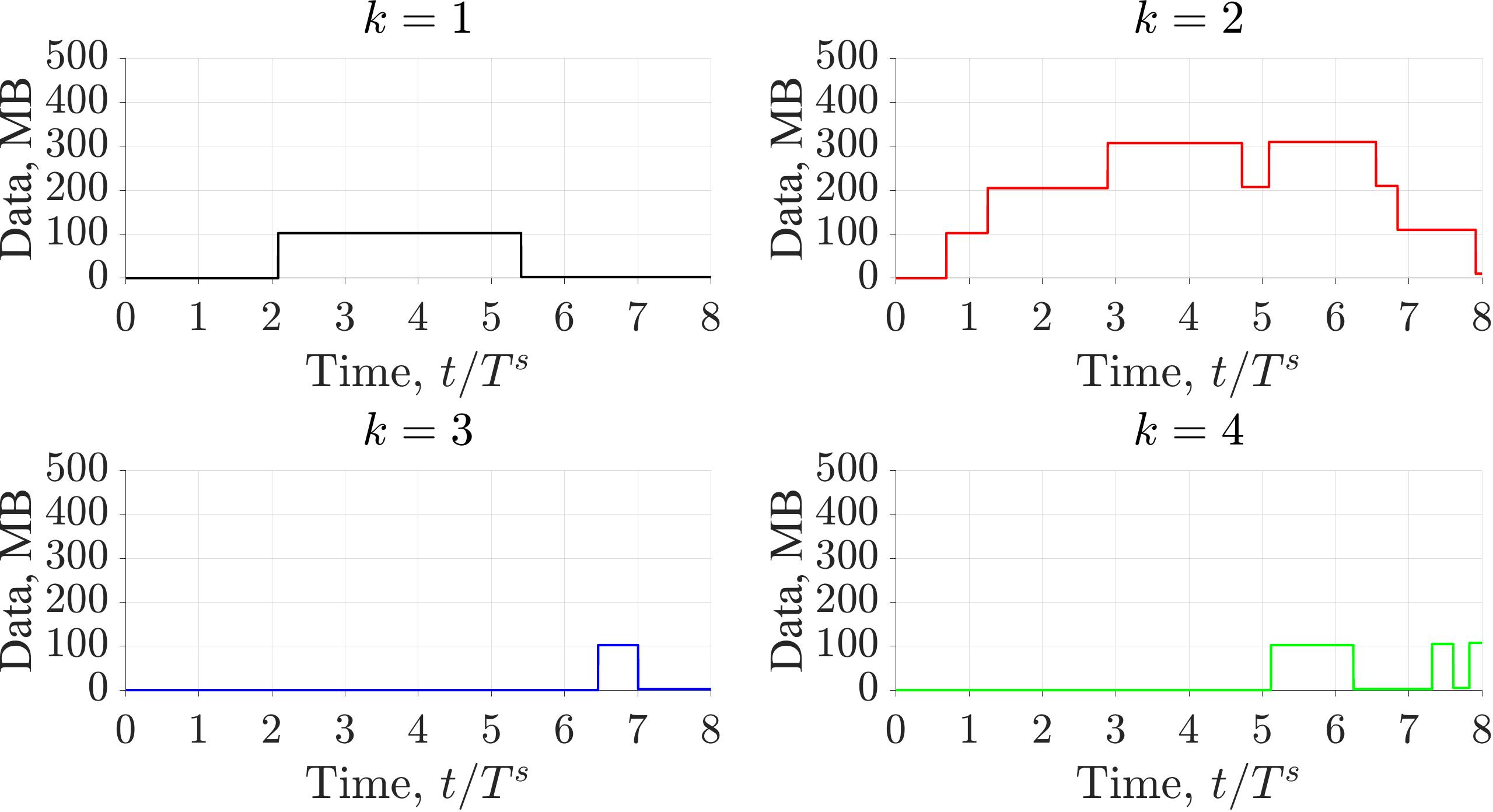}
        \caption{EOSSP data tracking.}
        \label{fig:EOSSP_Data}
    \end{subfigure}
    \begin{subfigure}[h]{0.49\textwidth}
        \centering
        \includegraphics[width = \textwidth]{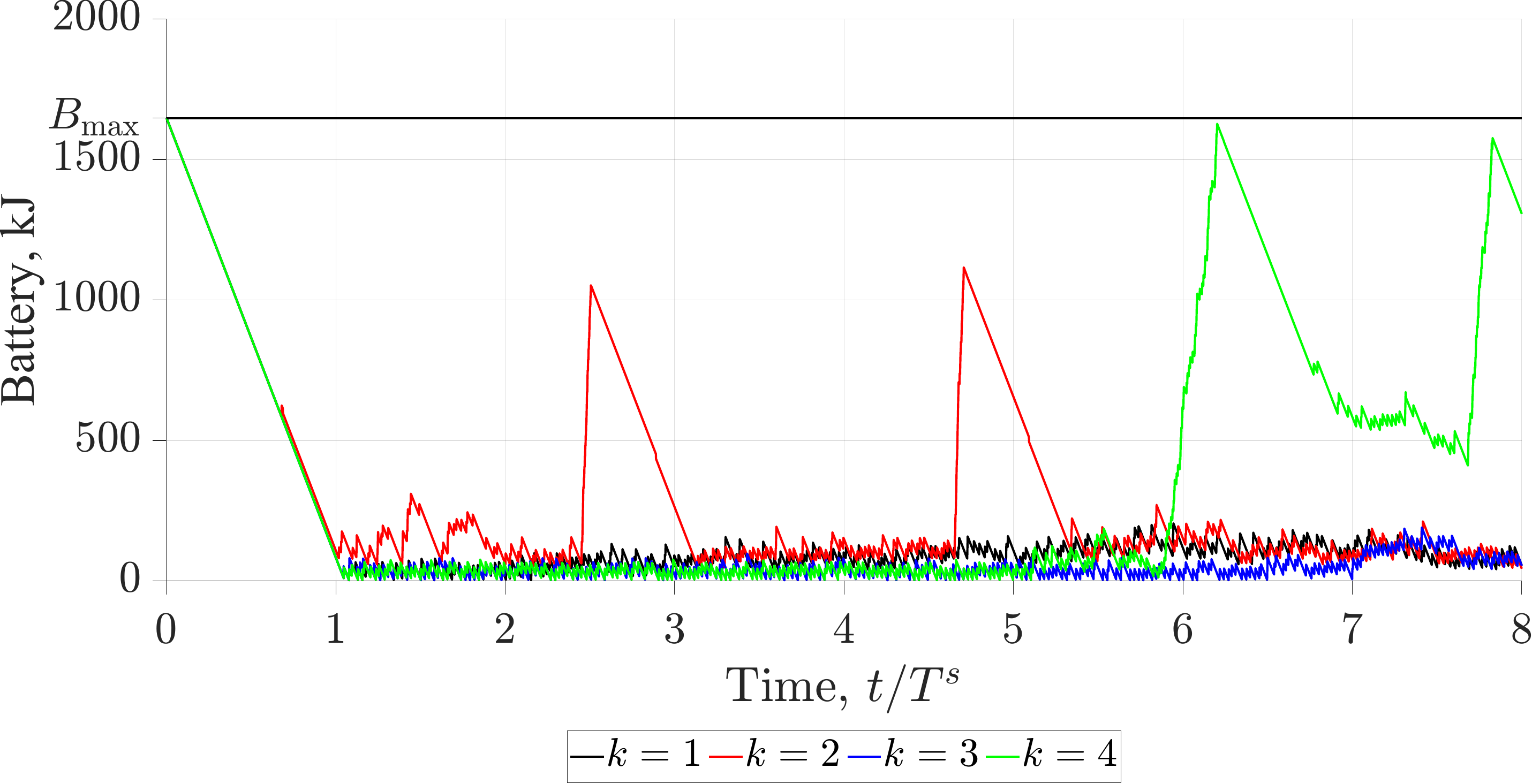}
        \caption{EOSSP battery tracking.}
        \label{fig:EOSSP_BATT}
    \end{subfigure}
    \caption{EOSSP schedule.}
    \label{fig:EOSSP_Schedule}
\end{figure}

\begin{figure}[!ht]
    \centering
    \begin{subfigure}[h]{0.49\textwidth}
        \centering
        \includegraphics[width = \textwidth]{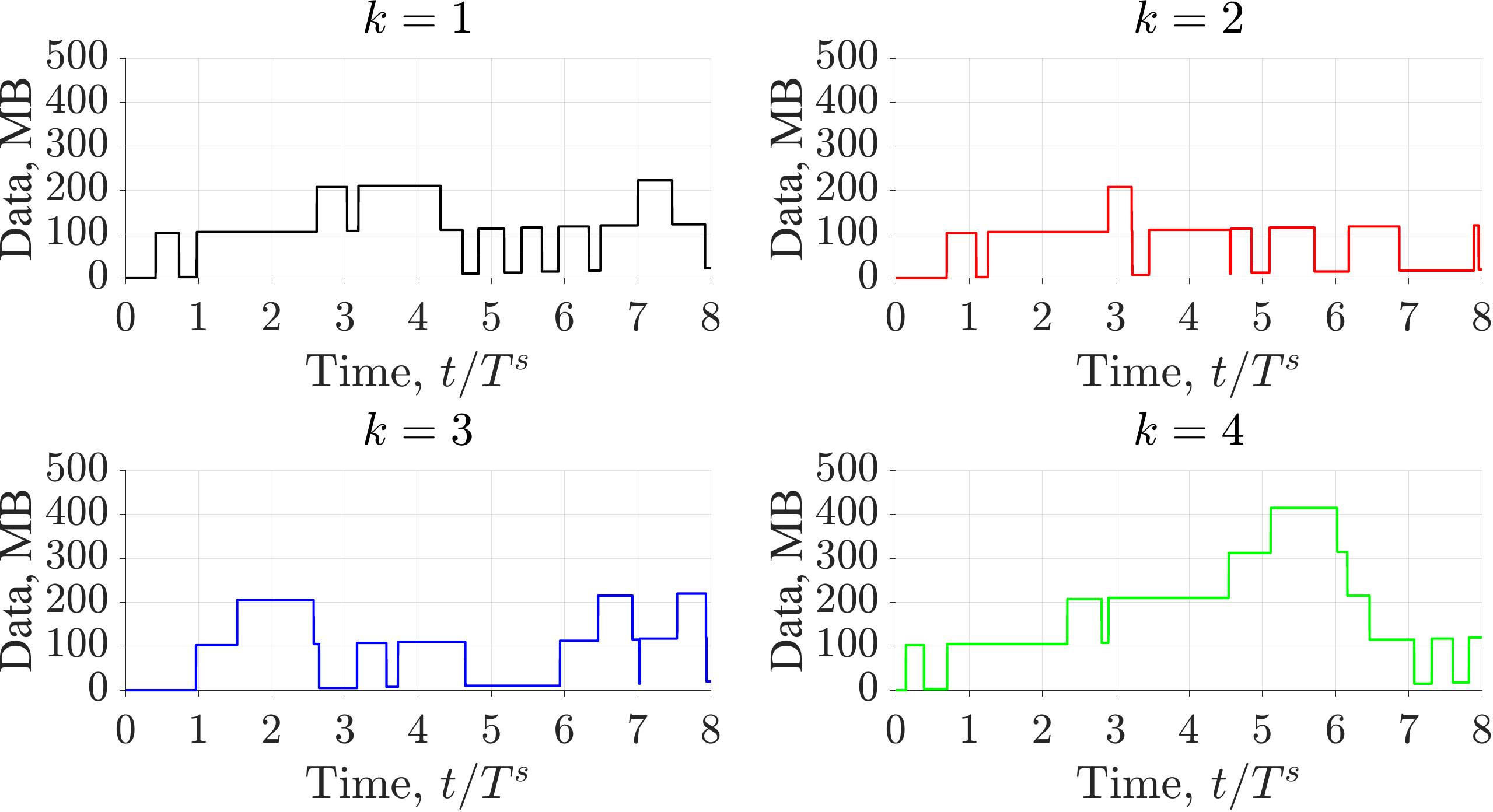}
        \caption{REOSSP-Exact data tracking}
        \label{fig:REOSSP_Data}
    \end{subfigure}
    \begin{subfigure}[h]{0.49\textwidth}
        \centering
        \includegraphics[width = \textwidth]{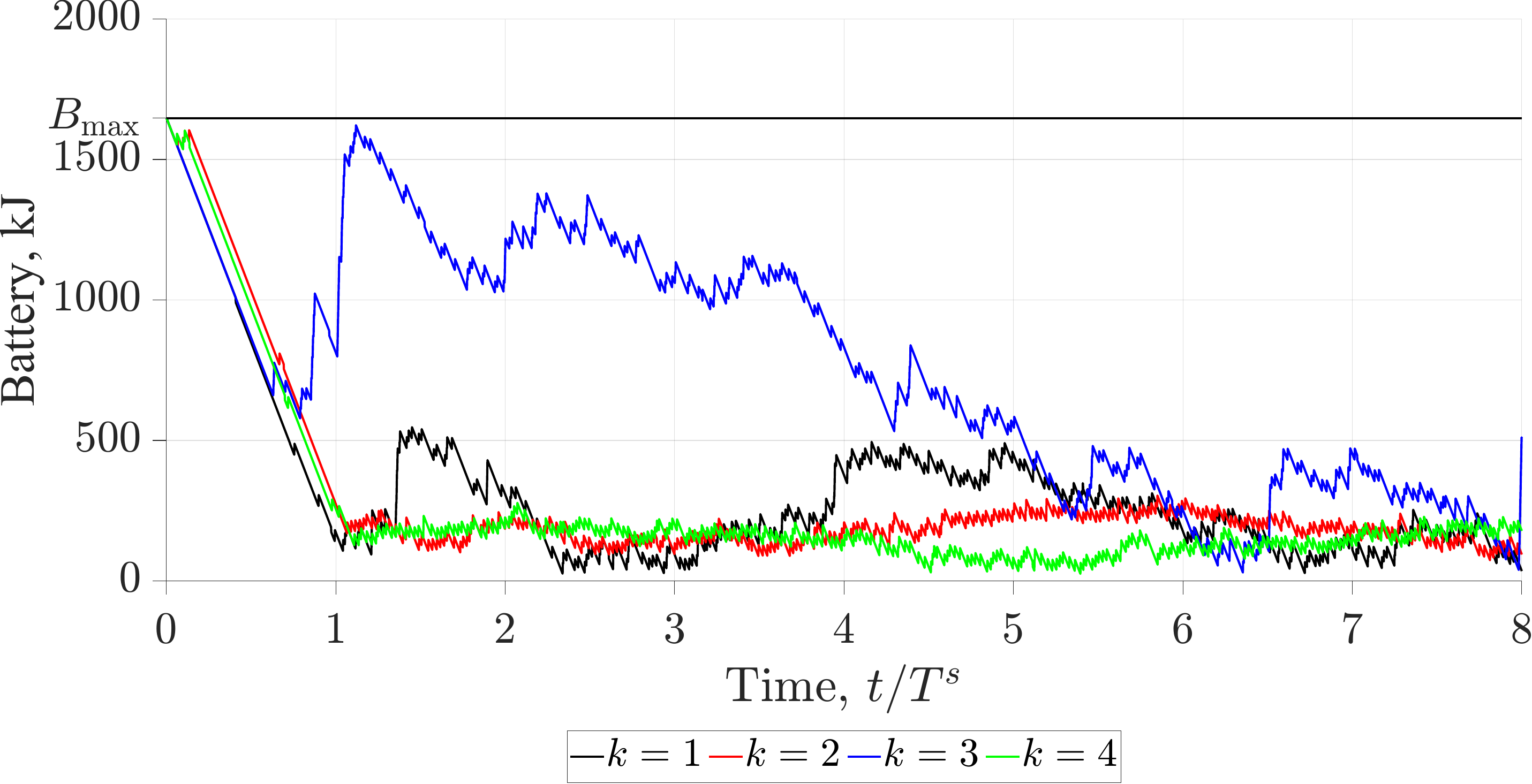}
        \caption{REOSSP-Exact battery tracking}
        \label{fig:REOSSP_BATT}
    \end{subfigure}
    \caption{REOSSP-Exact schedule.}
    \label{fig:REOSSP_Schedule}
\end{figure}

\begin{figure}[!ht]
    \centering
    \begin{subfigure}[h]{0.49\textwidth}
        \centering
        \includegraphics[width = \textwidth]{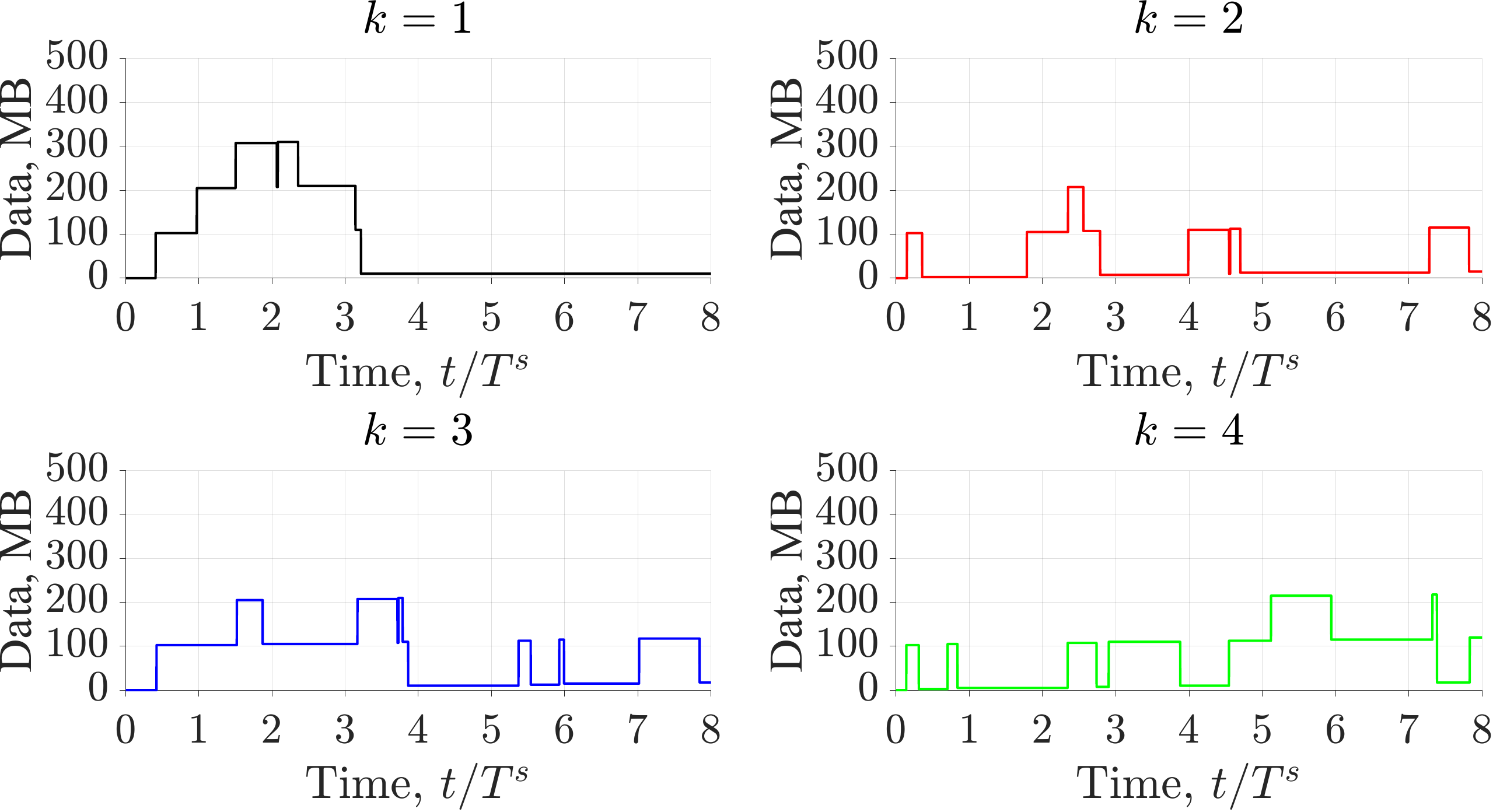}
        \caption{REOSSP-RHP data tracking.}
        \label{fig:RHP_Data}
    \end{subfigure}
    \begin{subfigure}[h]{0.49\textwidth}
        \centering
        \includegraphics[width = \textwidth]{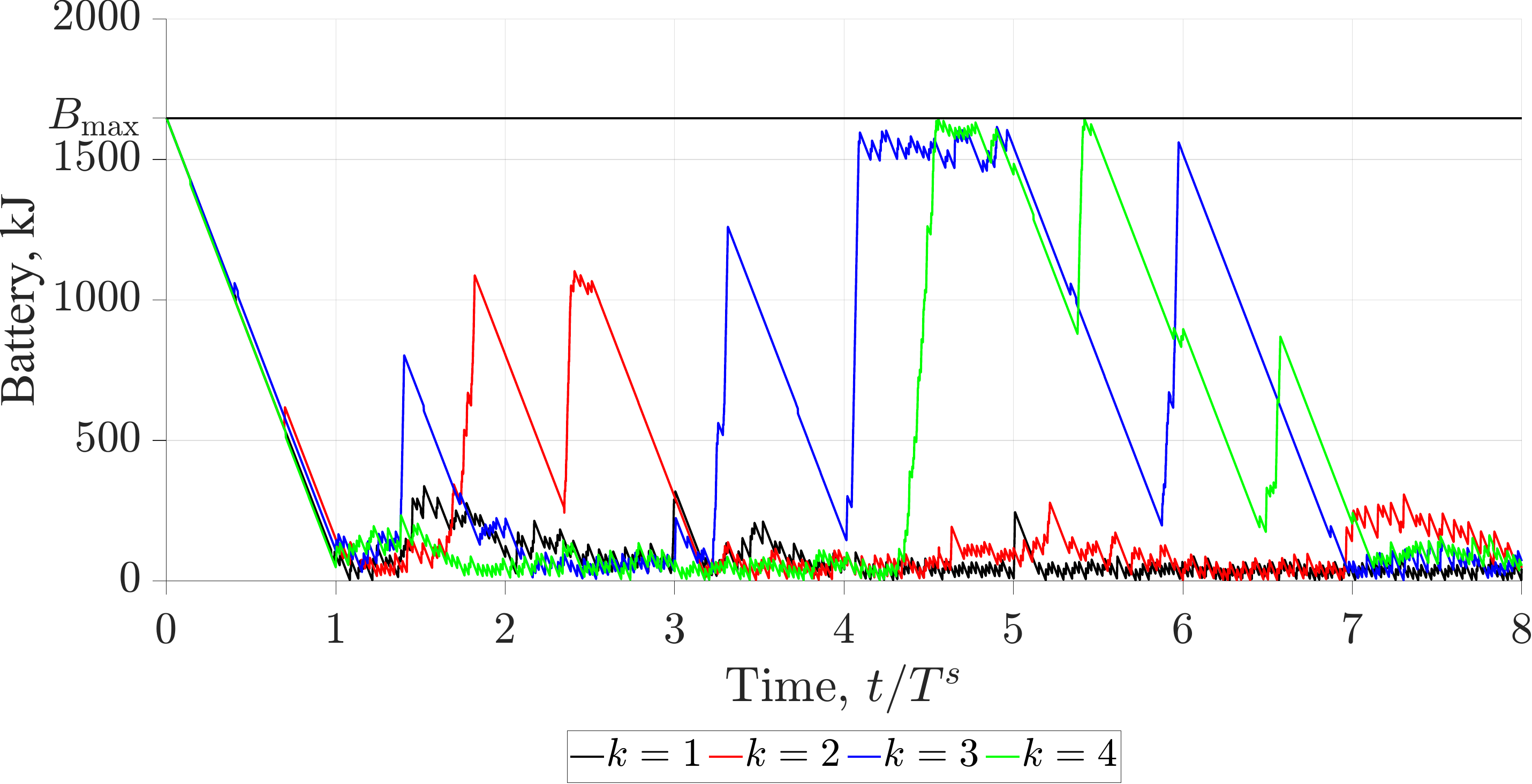}
        \caption{REOSSP-RHP battery tracking.}
        \label{fig:RHP_BATT}
    \end{subfigure}
    \caption{REOSSP-RHP schedule.}
    \label{fig:RHP_Schedule}
\end{figure}

The resultant schedules of the \EOSSP, \REOSSPExact, and \RHP are shown in Tables~\ref{tab:EOSSP_obs-down-charge}--\ref{tab:RHP_obs-down-charge}, respectively. The tables depict the high-level schedule broken down by stage (in terms of the time in which a stage would occur for the \EOSSP), including the number of observations, downlink occurrences, and charging occurrences, as well as the various orbital elements (inclination $i$, RAAN $\Omega$, and argument of latitude $u$, with changed elements highlighted in bold in the arrival stage) and the cost of orbital maneuvers for the \REOSSPExact and \RHP. It should be noted that Tables~\ref{tab:REOSSP_obs-down-charge} and~\ref{tab:RHP_obs-down-charge} report the values of $i$, $\Omega$, and $u$ belonging to a specific orbital slot as defined at the start of the schedule, although propagation with SGP4 accounts for perturbations throughout the schedule duration, including in the generation of visibility and orbital maneuver cost parameters. Furthermore, the data and battery storage capacities over time are shown for the \EOSSP, \REOSSPExact, and \RHP in Figs.~\ref{fig:EOSSP_Schedule}--\ref{fig:RHP_Schedule}, respectively. The objective function values of the \EOSSP, \REOSSPExact, and \RHP are $z_{E} = 25$, $z_{R} = 97$, and $\sum_{s = 1}^{S} z_{\text{RHP}} = 73$, respectively, while the figure of merit values for the \EOSSP, \REOSSPExact, and \RHP are $0.80$, $3.20$, and \SI{2.40}{GB}, respectively. It should be noted that despite the relatively low amount of data gathered by the constellation within the schedule horizon, the schedules only concern a single target (Hurricane Sandy) rather than any other passive monitoring occurring. Additionally, the runtime of the \EOSSP, \REOSSPExact, and \RHP are $0.24$, $54.60$, and \SI{6.65}{min}, respectively. As such, the \REOSSPExact and \RHP improved upon the \EOSSP in regards to the objective function value by \SI{288.00}{\%} and \SI{192.00}{\%} respectively, and in regards to the figure of merit by \SI{300.00}{\%} and \SI{200.00}{\%}, respectively. Separately, the \RHP performed \SI{24.74}{\%} and \SI{25.00}{\%} worse than the \REOSSPExact regarding the objective function value and figure of merit, respectively. Meanwhile, the runtime of the \RHP was only \SI{12.18}{\%} the runtime of the \REOSSPExact, substantially improving computation time at a slight performance cost. 

\begin{table}[!ht]
    \centering
    \caption{EOSSP observation, downlink, and charge occurrence by stage}
    \resizebox{0.75\textwidth}{!}{
    \begin{tabular}{ l l r r r r r r r r }
         \hline \hline
         Satellite & State & $S = 1$ & $S = 2$ & $S = 3$ & $S = 4$ & $S = 5$ & $S = 6$ & $S = 7$ & $S = 8$ \\
         \hline 
         & Observations & 0 & 0  & 1  & 0  & 0  & 0  & 0  & 0  \\
1        & Downlinks    & 0 & 0  & 0  & 0  & 0  & 1  & 0  & 0  \\
         & Charging     & 0 & 36 & 40 & 38 & 38 & 38 & 37 & 37 \\
         \hline
         & Observations & 1 & 1  & 1  & 0  & 0  & 1  & 0  & 0  \\
2        & Downlinks    & 0 & 0  & 0  & 0  & 1  & 0  & 2  & 1  \\
         & Charging     & 1 & 37 & 43 & 35 & 50 & 26 & 36 & 37 \\
         \hline
         & Observations & 0 & 0  & 0  & 0  & 0  & 0  & 1  & 0  \\
3        & Downlinks    & 0 & 0  & 0  & 0  & 0  & 0  & 0  & 1  \\
         & Charging     & 0 & 37 & 37 & 39 & 37 & 37 & 39 & 38 \\
         \hline
         & Observations & 0 & 0  & 0  & 0  & 0  & 1  & 0  & 2  \\
4        & Downlinks    & 0 & 0  & 0  & 0  & 0  & 0  & 1  & 1  \\
         & Charging     & 0 & 37 & 38 & 37 & 38 & 52 & 37 & 56 \\
         \hline \hline 
    \end{tabular}
    }
    \label{tab:EOSSP_obs-down-charge}
\end{table}

The optimal schedule of the \EOSSP contains only at most four observations per satellite, with two satellites only obtaining a single observation, which results in a low quantity of data and as such an objective function value of $z_{E} = 25$ and a figure of merit value of $Z_{E} = \SI{0.80}{GB}$. Furthermore, since the objective is to balance the observation and downlink of data, satellite 4 achieves two observations in stage eight while only downlinking one of these observations, resulting in the satellite retaining some data that may be downlinked at a time later than the schedule horizon established. The extra data retained by satellites 1 through 3 is a result of the previously mentioned discrepancy between the values $D_{\text{obs}}$ and $D_{\text{comm}}$, a difference of \SI{2.50}{MB} extra in $D_{\text{obs}}$. The extra data remaining onboard each satellite is reflected in Fig.~\ref{fig:EOSSP_data} where satellite 4 retains \SI{107.50}{MB}, satellite 2 retains \SI{10.00}{MB}, and satellites 1 and 3 retain \SI{2.50}{MB}. In addition, the charging of each satellite is dominated by the battery constraints due to the small number of tasks performed in each stage and the large amount of power obtained via charging, as shown in Fig.~\ref{fig:EOSSP_BATT}. However, the large peaks shown are not determined by an eclipse of sunlight, as any eclipse only lasts at most \SI{35}{min}, slightly more than one-third of an orbital period. The reason for the large peaks shown is the large variety of feasible charging opportunities due to the low occurrence of an eclipse of sunlight, resulting in some charging schedules including the large peaks (such as satellites 2 and 4) while others do not (such as satellites 1 and 3). 

\begin{table}[!ht]
    \centering
    \caption{REOSSP-Exact observation, downlink, and charge occurrence by stage.}
    \resizebox{0.75\textwidth}{!}{
    \begin{tabular}{ l l r r r r r r r r }
         \hline \hline
         Satellite & State & $S = 1$ & $S = 2$ & $S = 3$ & $S = 4$ & $S = 5$ & $S = 6$ & $S = 7$ & $S = 8$ \\
         \hline 
         & $i$, deg      & 98.18  & 98.18  & 98.18  & 98.18  & 98.18  & 100.33 & 100.33 & 100.33 \\
         & $\Omega$, deg & 0.00   & 0.00   & 0.00   & 0.00   & 0.00   & 0.00   & 0.00   & 0.00   \\
         & $u$, deg      & 240    & 240    & 240    & 192    & 192    & 48     & 312    & 312    \\
1        & Observations  & 2      & 0      & 1      & 1      & 1      & 2      & 2      & 0      \\
         & Donwlinks     & 1      & 0      & 0      & 1      & 2      & 2      & 1      & 2      \\
         & Charging      & 3      & 41     & 33     & 45     & 40     & 32     & 37     & 36     \\
         & Cost, m/s     & 189.94 & 0.00   & 0.00   & 66.13  & 0.00   & 365.01 & 101.67 & 0.00   \\
         \hline 
         & $i$, deg      & 98.18  & 98.18  & 98.18  & 98.18  & 100.33 & 100.33 & 100.33 & 100.33 \\
         & $\Omega$, deg & 90.00  & 90.00  & 90.00  & 90.00  & 90.00  & 90.00  & 90.00  & 90.00  \\
         & $u$, deg      & 336    & 336    & 336    & 336    & 72     & 336    & 240    & 168    \\
2        & Observations  & 1      & 1      & 1      & 1      & 1      & 1      & 1      & 1      \\
         & Donwlinks     & 0      & 1      & 0      & 2      & 2      & 1      & 1      & 1      \\
         & Charging      & 5      & 37     & 36     & 40     & 39     & 39     & 36     & 36     \\
         & Cost, m/s     & 145.96 & 0.00   & 0.00   & 0.00   & 347.48 & 101.67 & 101.67 & 16.37  \\
         \hline 
         & $i$, deg      & 98.18  & 98.18  & 98.18  & 98.18  & 98.18  & 98.18  & 98.18  & 98.18  \\
         & $\Omega$, deg & 177.83 & 177.83 & 177.83 & 177.83 & 177.83 & 177.83 & 177.83 & 177.83 \\
         & $u$, deg      & 144    & 144    & 144    & 144    & 72     & 72     & 0      & 288    \\
3        & Observations  & 1      & 1      & 0      & 2      & 0      & 1      & 1      & 2      \\
         & Donwlinks     & 0      & 0      & 2      & 1      & 1      & 0      & 1      & 3      \\
         & Charging      & 18     & 48     & 34     & 33     & 31     & 31     & 43     & 40     \\
         & Cost, m/s     & 461.58 & 0.00   & 0.00   & 0.00   & 16.37  & 0.00   & 16.37  & 16.37  \\
         \hline 
         & $i$, deg      & 98.18  & 98.18  & 98.18  & 98.18  & 98.18  & 98.18  & 98.18  & 98.18  \\
         & $\Omega$, deg & 267.83 & 267.83 & 267.83 & 267.83 & 267.83 & 267.83 & 267.83 & 267.83 \\
         & $u$, deg      & 96     & 96     & 96     & 96     & 96     & 24     & 24     & 24     \\
4        & Observations  & 2      & 0      & 2      & 0      & 1      & 1      & 0      & 2      \\
         & Donwlinks     & 1      & 0      & 1      & 0      & 0      & 0      & 3      & 2      \\
         & Charging      & 5      & 37     & 38     & 36     & 37     & 39     & 39     & 39     \\
         & Cost, m/s     & 301.57 & 0.00   & 0.00   & 0.00   & 0.00   & 16.37  & 0.00   & 0.00   \\
         \hline \hline 
    \end{tabular}
    }
    \label{tab:REOSSP_obs-down-charge}
\end{table}

The optimal schedule of the \REOSSPExact performs significantly better than the \EOSSP, by \SI{288.00}{\%} in terms of the objective function value and \SI{300.00}{\%} in terms of the figure of merit, including at least eight observations per satellite. The \REOSSPExact also suffers from an observation with no subsequent downlink in stage eight by satellite 4. As such, the extra data retained by each satellite is slightly higher than that retained by the satellites in the \EOSSP as a result of the higher quantity of observations. Specifically reflected in Fig.~\ref{fig:REOSSP_Data}, satellite 4 retains \SI{120.00}{MB}, satellite 1 retains \SI{22.50}{MB}, and satellites 2 and 3 retain \SI{20.00}{MB}. Additionally, satellite 4 reaches the highest data storage value of \SI{415.00}{MB} within the fifth stage, wherein satellite 4 downlinks all data before the final stage (stage eight) at which point two more observations are conducted with only a single downlink, which accounts for the extra data retained. The charging schedule, shown in Fig.~\ref{fig:REOSSP_BATT}, of the \REOSSPExact may appear vastly different from the previous charging schedule; however, this variation is a result of the large variety of feasible charging solutions as the longest eclipse for this orbital maneuver pattern is similar to that of the \EOSSP. Finally, the overall cost of the optimal \REOSSPExact schedule consumes \SI{2.26}{km/s}, \SI{75.48}{\%} of the total budget provided, for a total of 15 transfers (of 32 possibilities for transfers), of which four are a plane change (two of which are inclination raising and two of which are RAAN lowering). Each satellite, listed in order 1 through 4, consumed \SI{722.75}{m/s}, \SI{713.15}{m/s}, \SI{510.69}{m/s}, and \SI{317.93}{m/s}, which is \SI{96.37}{\%}, \SI{95.09}{\%}, \SI{68.09}{\%}, and \SI{42.39}{\%} of each individual budget, respectively. 

\begin{table}[!ht]
    \centering
    \caption{REOSSP-RHP observation, downlink, and charge occurrence by stage.}
    \resizebox{0.75\textwidth}{!}{
    \begin{tabular}{ l l r r r r r r r r }
         \hline \hline
         Satellite & State & $S = 1$ & $S = 2$ & $S = 3$ & $S = 4$ & $S = 5$ & $S = 6$ & $S = 7$ & $S = 8$ \\
         \hline 
         & $i$, deg      & 98.18  & 98.18  & 98.18  & 98.18  & 98.18  & 98.18  & 98.18  & 98.18  \\
         & $\Omega$, deg & 0.00   & 357.83 & 357.83 & 357.83 & 357.83 & 357.83 & 357.83 & 357.83 \\
         & $u$, deg      & 240    & 96     & 24     & 312    & 312    & 312    & 312    & 312    \\
1        & Observations  & 2      & 1      & 1      & 0      & 0      & 0      & 0      & 0      \\
         & Donwlinks     & 0      & 0      & 2      & 2      & 0      & 0      & 0      & 0      \\
         & Charging      & 1      & 39     & 42     & 32     & 37     & 39     & 37     & 38     \\
         & Cost, m/s     & 189.94 & 523.08 & 16.37  & 16.37  & 0.00   & 0.00   & 0.00   & 0.00   \\
         \hline 
         & $i$, deg      & 93.88  & 93.88  & 93.88  & 93.88  & 93.88  & 93.88  & 93.88  & 93.88  \\
         & $\Omega$, deg & 90.00  & 90.00  & 90.00  & 90.00  & 90.00  & 90.00  & 90.00  & 90.00  \\
         & $u$, deg      & 192    & 192    & 192    & 192    & 120    & 72     & 72     & 72     \\
2        & Observations  & 1      & 1      & 1      & 1      & 1      & 0      & 0      & 1      \\
         & Donwlinks     & 1      & 0      & 2      & 0      & 2      & 0      & 0      & 1      \\
         & Charging      & 2      & 54     & 26     & 33     & 39     & 36     & 42     & 34     \\
         & Cost, m/s     & 603.92 & 0.00   & 0.00   & 0.00   & 16.37  & 66.13  & 0.00   & 0.00   \\
         \hline 
         & $i$, deg      & 98.18  & 96.03  & 96.03  & 96.03  & 96.03  & 96.03  & 96.03  & 96.03  \\
         & $\Omega$, deg & 180.00 & 180.00 & 180.00 & 180.00 & 180.00 & 180.00 & 180.00 & 180.00 \\
         & $u$, deg      & 336    & 168    & 120    & 120    & 120    & 120    & 48     & 336    \\
3        & Observations  & 1      & 1      & 0      & 2      & 0      & 2      & 0      & 1      \\
         & Donwlinks     & 0      & 1      & 0      & 3      & 0      & 2      & 0      & 1      \\
         & Charging      & 1      & 41     & 34     & 41     & 71     & 38     & 2      & 39     \\
         & Cost, m/s     & 145.96 & 466.27 & 66.13  & 0.00   & 0.00   & 0.00   & 16.37  & 16.37  \\
         \hline 
         & $i$, deg      & 98.18  & 98.18  & 98.18  & 98.18  & 98.18  & 98.18  & 98.18  & 98.18  \\
         & $\Omega$, deg & 270.00 & 270.00 & 270.00 & 270.00 & 270.00 & 270.00 & 270.00 & 272.17 \\
         & $u$, deg      & 72     & 72     & 72     & 72     & 72     & 0      & 0      & 336    \\
4        & Observations  & 2      & 0      & 2      & 0      & 1      & 1      & 0      & 2      \\
         & Donwlinks     & 2      & 0      & 1      & 1      & 0      & 1      & 0      & 2      \\
         & Charging      & 1      & 37     & 39     & 38     & 72     & 24     & 22     & 34     \\
         & Cost, m/s     & 441.25 & 0.00   & 0.00   & 0.00   & 0.00   & 16.37  & 0.00   & 291.96 \\
         \hline \hline 
    \end{tabular}
    }
    \label{tab:RHP_obs-down-charge}
\end{table}

The optimal schedule of the \RHP also performs better than the \EOSSP, by \SI{192.00}{\%} in terms of the objective function value and \SI{200.00}{\%} in terms of the figure of merit, but slightly worse than the \REOSSPExact, with at least four observations per satellite. The \RHP also suffers from an observation with no subsequent downlink in stage eight, indicating that this phenomenon is a result of the visibility profiles of the target and ground stations, wherein satellite four has visibility of a target near the end of the time horizon with no visibility of ground stations afterward. Since the overall number of observations conducted by the \RHP is less than that of the \REOSSPExact, the extra data retained by each satellite is $120.00$, $17.50$, $15.00$, and \SI{10.00}{MB} for satellites 4 through 1, respectively, as reflected in Fig.~\ref{fig:RHP_Data}. An additional phenomenon occurs in Fig.~\ref{fig:RHP_Data}, wherein data that is observed by each satellite is downlinked to a ground station much sooner than in the \REOSSPExact and \EOSSP, a result of the objective to downlink and observe being optimized for each group of two stages rather than over the accumulation of the entire time horizon. A similar phenomenon occurs in Fig.~\ref{fig:RHP_BATT}, wherein the charging of each satellite is maintained within each subproblem, leading to a slightly more sporadic charging schedule when appended. Finally, the overall cost of the optimal \RHP schedule is substantially more expensive than that of the \REOSSPExact, a result of the rather small $L=1$ stage lookahead incentivizing effective yet expensive maneuvers to occur at the start of the schedule. This is displayed well in Table~\ref{tab:RHP_obs-down-charge}, where each satellite consumes on average \SI{79.01}{\%} of the maneuver budget within the first two maneuvers. Overall, the optimal \RHP schedule consumes \SI{2.89}{km/s}, \SI{96.43}{\%} of the total budget provided, for a total of 15 transfers (of 32 possibilities for transfers), of which four are a plane change (two of which are inclination lowering, one of which is RAAN lowering, and one of which is RAAN raising). Each satellite, listed in order 1 through 4, consumed \SI{745.76}{m/s}, \SI{686.42}{m/s}, \SI{711.10}{m/s}, and \SI{749.58}{m/s}, which is \SI{99.43}{\%}, \SI{91.52}{\%}, \SI{94.81}{\%}, and \SI{99.94}{\%} of each individual budget, respectively. 

Overall, the highly adaptable nature of the \REOSSPExact and \RHP enabled through reconfigurability, a diverse set of orbital slots including plane and phase changes, and the reduction in computation time via the \RHP algorithm allows for performance much higher than that of the \EOSSP. This effectiveness is shown through the performance of significantly more observations and downlink occurrences by each satellite relative to the \EOSSP, while the \RHP performs more expensive orbital maneuvers while completing computations much faster. Additionally, all three scheduling algorithms return a resultant schedule that obeys all constraints, including power consumption, charging, data management, and orbital maneuver budget limitations. 

\section{Conclusions} \label{sec:Conclusion}

This paper developed a mathematical formulation to integrate constellation reconfigurability into the EOSSP through the use of MILP and previous research on the MCRP \cite{Pearl2023Comparing,Pearl2024Benchmarking,Lee2022,Lee2023,lee2024deterministic}. The \REOSSP developed herein optimizes the schedule of a cooperative satellite constellation capable of orbital maneuvering such that the number of observations and subsequent data downlink occurrences is maximized, with a higher influence from the number of downlink occurrences. Additionally, the \REOSSP considers available data storage, battery storage, and a propellant budget for orbital maneuvers, while scheduling the tasks to observe targets, downlink observed data, charge onboard batteries, and perform constellation reconfiguration at predetermined stages. Furthermore, the \REOSSP is solved via the \RHP, which significantly improves computational runtime compared to the \REOSSP but at the cost of yielding suboptimal solutions. The novel MILP formulation of the \REOSSP and the \RHP allows flexible implementation, including the ability to include trajectory optimization for orbital maneuvers and optimal solutions through commercial solvers. In conjunction, a baseline \EOSSP is developed and utilized for the comparison. 

The computational experiments conducted in Sec.~\ref{sec:Experiments} demonstrate the capability of the \REOSSPExact to outperform the \EOSSP in a variety of scenarios, confirming the hypothesis set forth by previous comparative research \cite{Pearl2023Comparing,Pearl2024Benchmarking,lee2024deterministic}. Within the random instances, meant to demonstrate the flexibility of the framework under various parameters, the average percent improvement of the \REOSSPExact and \RHP over the \EOSSP is \SI{101.80}{\%} and \SI{78.06}{\%} in terms of the objective function value, respectively. When applied to Hurricane Sandy the \REOSSPExact and \RHP exceed the performance of the \EOSSP by \SI{288.00}{\%} and \SI{192.00}{\%} in terms of the objective function value, respectively, and \SI{300.00}{\%} and \SI{200.00}{\%} in terms of the figure of merit, respectively. While the \RHP does not obtain as high-quality results as the \REOSSPExact, in the event that the problem scale is sufficiently large so as to cause the \REOSSPExact to be unable to converge to a feasible solution, the \RHP can obtain an optimal solution to each subproblem while still significantly outperforming the \EOSSP. Both the random instances and the case study depict the performance increase provided by the \REOSSPExact and \RHP, while simultaneously obeying all operational and physical constraints. Overall, the implementation of constellation reconfigurability provides a greater throughput of observed data than the standard \EOSSP. 

Future research into the concept of constellation reconfigurability implemented in a scheduling format may take many forms. Firstly, the scope of this paper is restricted to high-thrust orbital maneuver trajectory optimization, whereas low-thrust maneuvers may allow for a lower operational cost of such a reconfigurable constellation. As such, the implementation of various low-thrust trajectory optimization theories, such as those in Refs.~\cite{Sims-Flanagan2000,Jackson2020QLaw}, may result in an increased level of performance above the baseline with less expensive orbital maneuvers. Secondly, all scheduling problems formulated in this paper assume deterministic characteristics, while many natural disasters or other observable planetary phenomena are often stochastic. Therefore, the implementation of stochastic programming similar to that in Ref.~\cite{Hoskins2024} to the scheduling problems developed in this paper may more accurately represent the stochastic nature of natural disaster targets. Thirdly, the work in this paper focuses heavily on the \REOSSP and constellation reconfigurability as a whole, while previous work in Refs.~\cite{Pearl2023Comparing,Pearl2024Benchmarking} additionally considers satellite agility, thus contributing to the potential of a similar AEOSSP formulation to be developed and used for further benchmarking of satellite agility against constellation reconfigurability. Finally, the further comparison of a small constellation of maneuverable satellites against a much larger fixed constellation of either nadir-directional or agile satellites may additionally prove insightful.

\section*{Appendix A: Proof of Formulation Equivalence with No Orbital Maneuvers}

The \EOSSP and \REOSSP are based on the same principles, such that the only differences lie within the orbital maneuverability, wherein the \REOSSP contributes the properties of the currently occupied orbital slots. Such properties include the visibility of task priorities and the battery required to perform orbital maneuvers, $B_{\text{recon}}$.

\textit{Theorem $1$:} Let $c^k_{\max} = 0, \forall k \in \mathcal{K}$ in the \REOSSP, then $z_{E} = z_{R}$.

\textit{Proof:} Suppose that for the \EOSSP and \REOSSP, all constant parameters such as $\Delta t$ or $K$, are identical. Note that each stage aligns sequentially such that $t \in \mathcal{T}^s$ of the \REOSSP corresponds to $t + (s-1)T^s \in \mathcal{T}$ of the \EOSSP, as reflected in Fig.~\ref{fig:recon_demonstration}. This sequential alignment also applies to the visibility parameters $V^{sk}_{tjp}$, $G^{sk}_{tjg}$, and $H^{sk}_{tj}$. Additionally, the initial orbital positions of each satellite (those applied to the \EOSSP) apply to a single orbital slot option index $j$. Finally, the number of stages $S$ and the number of orbital slots $J^{sk}, \forall s \in \mathcal{S} \setminus \{0\}, \forall k \in \mathcal{K}$ may remain generic. 

Firstly, constraints~\eqref{REOSSP:flow_cost} restrict the orbital maneuvers of satellites to those with a cost $c^{sk}_{ij}$ whose accumulation does not exceed $c^k_{\max}$. Since $c^k_{\max} = 0, \forall k \in \mathcal{K}$, only the maneuvers from orbital slot $i \in \mathcal{J}^{s-1, k}$ to orbital slot $j \in \mathcal{J}^{sk}$ where $c^{sk}_{ij} = 0$, which only occurs if $i = j, \forall k \in \mathcal{K}, \forall s \in \mathcal{S} \setminus \{0\}$ by definition, are feasible. Furthermore, at stage $s=1$, the orbital slot $i \in \mathcal{J}^{0k}$ is a singleton set of the initial orbital positions for each satellite (those applied to the \EOSSP). As such, the satellites of the \EOSSP and \REOSSP are located in the same orbital slot throughout the entire time horizon. Finally, the assignment $i \in \mathcal{J}^{0k}, \forall k \in \mathcal{K}$ and $j = i, j \in \mathcal{J}^{sk}, \forall k \in \mathcal{K}, \forall s \in \mathcal{S} \setminus \{0\}$ obeys constraints~\eqref{REOSSP:flow_from_initial_conditions} and~\eqref{REOSSP:flow_from0<s<S} since only orbital slot $i$ is occupied in all stages $s \in \mathcal{S}$ such that $x^{sk}_{ii} = 1, \forall k \in \mathcal{K}, \forall s \in \mathcal{S} \setminus \{0\}$ and $x^{sk}_{ij} = 0, \forall k \in \mathcal{K}, \forall s \in \mathcal{S} \setminus \{0\}, j \neq i$. Since each satellite cannot perform an orbital maneuver, it can be assumed that $B_{\text{recon}}=0$, thus eliminating the possible difference in battery power consumption.

Secondly, constraints~(\ref{REOSSP:target_visibility}--\ref{REOSSP:sun_visibility}) assign the visibility of targets, ground stations, and the sun, respectively. Since $i \in \mathcal{J}^{0k}, \forall k \in \mathcal{K}$ is the orbital position applied to the \EOSSP, $i = j, \forall k \in \mathcal{K}, \forall s \in \mathcal{S} \setminus \{0\}$, and the visibility parameters align sequentially, the visibility of targets, ground stations, and the sun are identically assigned to the \EOSSP and \REOSSP. Furthermore, the VTWs for the tasks of target observation, data downlink to ground stations, and solar charging are also identically assigned. Additionally, since objective function~\eqref{EOSSP:obj} of the \EOSSP and objective function~\eqref{REOSSP:obj} of the \REOSSP are functionally equivalent, and the VTWs of each task are identical, the assignment of tasks is subsequently identical. In a similar manner to the assignment of orbital slots, constraints~\eqref{REOSSP:obvs-down_overlap} are inherently obeyed as the identical task assignment obeys constraints~\eqref{EOSSP:obvs-down_overlap}.

Thirdly, constraints~\eqref{REOSSP:Data} track and restrict the data storage of each satellite to within the bounds $D_{\min}$ and $D_{\max}$ dependent upon the occurrence of target observation and data downlink to ground stations. Since the values of $D_{\text{comm}}$ and $D_{\text{obs}}$ are constant, each stage is sequentially aligned, and task occurrence is identical between the \EOSSP and \REOSSP, constraints~\eqref{REOSSP:d<max} and~\eqref{REOSSP:d>0} are inherently identical to constraints~\eqref{EOSSP:d<max} and~\eqref{EOSSP:d>0}, respectively, and are thusly obeyed. Additionally, constraints~\eqref{REOSSP:d-track_not_Ts} and~\eqref{REOSSP:d-track_Ts} contribute the current time step task performance to the subsequent time step data storage value, wherein constraints~\eqref{REOSSP:d-track_not_Ts} apply within each stage and constraints~\eqref{REOSSP:d-track_Ts} apply between each stage. As such, constraints~\eqref{REOSSP:d-track_not_Ts} and~\eqref{REOSSP:d-track_Ts} operate identically to constraints~\eqref{EOSSP:d-track}.

Fourthly, constraints~\eqref{REOSSP:Battery_Track} and~\eqref{REOSSP:Battery_Constrain} track and restrict the battery storage, respectively, of each satellite to within the bounds $B_{\min}$ and $B_{\max}$ dependent upon task performance; the consideration of orbital maneuvers may be removed through the assumption that $B_{\text{recon}} = 0$. Since the values of $B_{\text{comm}}$, $B_{\text{charge}}$, $B_{\text{obs}}$, and $B_{\text{time}}$ are constant, each stage is sequentially aligned, and task occurrence is identical, constraints~\eqref{REOSSP:b<max} are inherently identical to constraints~\eqref{EOSSP:b<max} and are also obeyed. As an aside, since $B_{\text{recon}}=0$, constraints~\eqref{REOSSP:b>0_Ts_s1} simplify to $B^k_{\max} \ge 0, \forall k \in \mathcal{K}$ which is implied, and constraints~\eqref{REOSSP:b>0_Ts_not_s1} simplify to $b^{sk}_{T^s} - \sum_{p \in \mathcal{P}}B_{\text{obs}}y^{sk}_{T^s p} - \sum_{g \in \mathcal{G}} B_{\text{comm}} q^{sk}_{T^s g} - B_{\text{time}} \ge 0, \forall s \in \mathcal{S} \setminus \{0, S\}, \forall k \in \mathcal{K}$. With these simplifications, constraints~\eqref{REOSSP:b>0_not_Ts} and~\eqref{REOSSP:b>0_Ts_not_s1} function similarly to constraints~\eqref{REOSSP:d-track_not_Ts} and~\eqref{REOSSP:d-track_Ts} in that both contribute the current time step task performance to the subsequent time step battery storage value within and between each stage, respectively. As such, constraints~\eqref{REOSSP:b>0_not_Ts} and~\eqref{REOSSP:b>0_Ts_not_s1} operate identically to constraints~\eqref{EOSSP:b>0}. Furthermore, constraints~\eqref{REOSSP:b-track_Ts_s1} simplify to $b^{1k}_1 = B^k_{\max}, \forall k \in \mathcal{K}$, an assumption of the \EOSSP, and constraints~\eqref{REOSSP:b-track_Ts_not_s1} simplify to $b^{s+1, k}_1 = b^{sk}_{T^s} + B_{\text{charge}} h^{sk}_{T^s} - \sum_{p \in \mathcal{P}} B_{\text{obs}} y^{sk}_{T^s p} - \sum_{g \in \mathcal{G}} B_{\text{comm}} q^{sk}_{T^s g} - B_{\text{time}}, \forall s \in \mathcal{S} \setminus \{0, S\}, \forall k \in \mathcal{K}$. Therefore, constraints~\eqref{REOSSP:b-track_not_Ts} and~\eqref{REOSSP:b-track_Ts_not_s1} function similarly to constraints~\eqref{REOSSP:d-track_not_Ts} and~\eqref{REOSSP:d-track_Ts}, and thus operate identically to constraints~\eqref{EOSSP:b-track}.

Finally, objective functions~\eqref{REOSSP:obj} and~\eqref{EOSSP:obj} compute the objective function values $z_{R}$ and $z_{E}$, respectively. Since task occurrence is identical and each stage is sequentially aligned, the values of $q^k_{tg}$ and $y^k_{tp}$ will be the same as the values of $q^{sk}_{tg}$ and $y^{sk}_{tp}$ for similar time steps $t + (s-1)T^s \in \mathcal{T}$ and $t \in \mathcal{T}^s$, respectively. As such, since the value of $C$ is constant, the objective function values $z_{E} = z_{R}$. Therefore, the \EOSSP and \REOSSP return identical schedules with identical objective function values under conditions wherein $c^k_{\max} = 0, \forall k \in \mathcal{K}$. \qedsymbol 

To further validate the above proof numerically, the case study in Sec.~\ref{subsec:case_study} is repeated with the parameter $c^k_{\max} = 0, \forall k \in \mathcal{K}$. The results follow that $z_{E} = z_{R} = 25$ since the satellites of the \REOSSP were forced to remain in the initial condition orbital slots for the entire duration of the schedule horizon. As such, the above proof is validated using realistic target data and the same parameters and parameter generation methods as the case study in Sec.~\ref{subsec:case_study}.

\section*{Appendix B: Additional Random Instance Results Statistics}

\begin{figure}[!ht]
    \centering
    \begin{subfigure}[h]{0.49\textwidth}
        \centering
        \includegraphics[width = \textwidth]{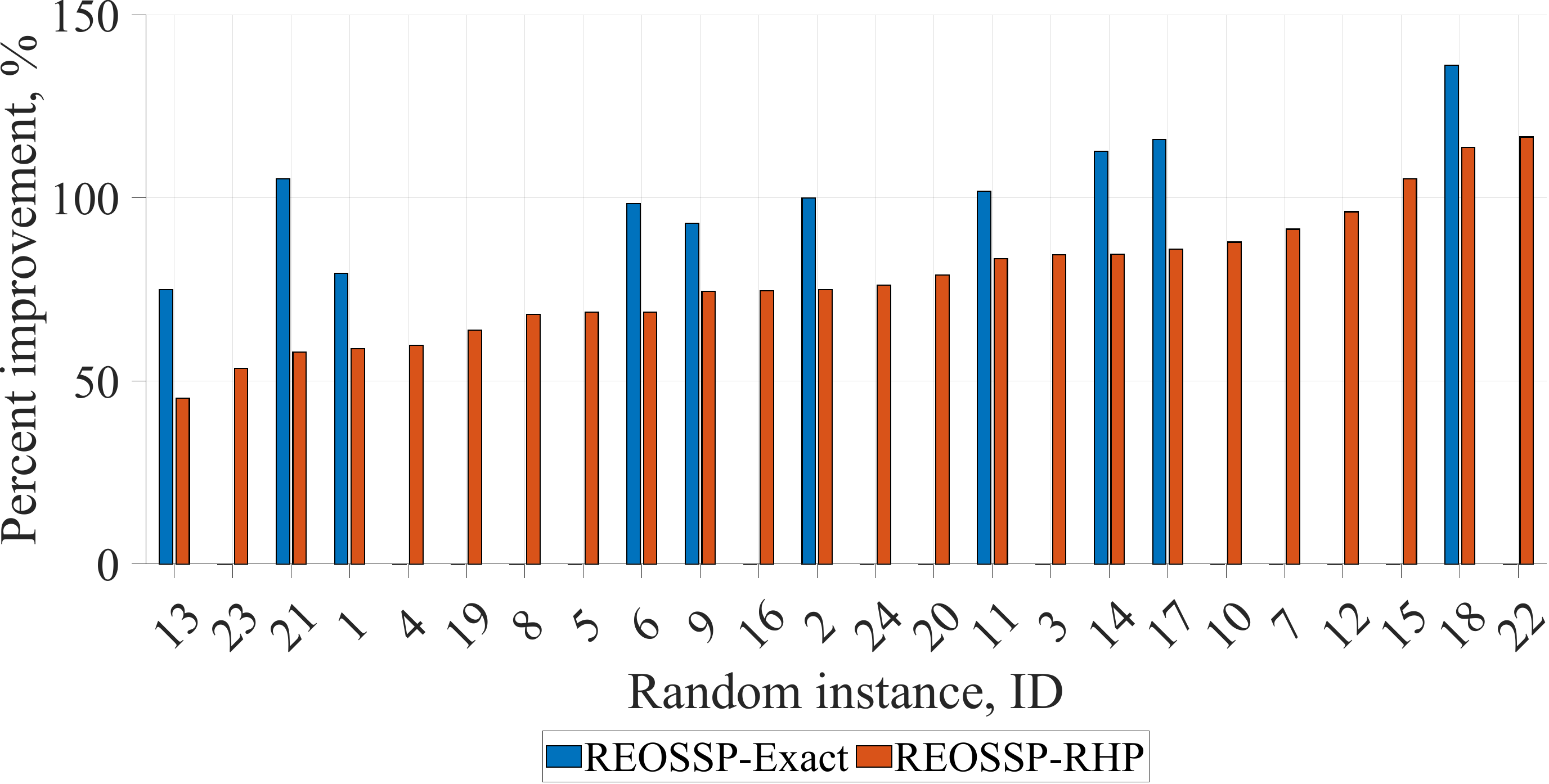}
        \caption{Percent improvement over the \EOSSP}
        \label{fig:random_percent}
    \end{subfigure}
    \hfill
    \begin{subfigure}[h]{0.49\textwidth}
        \centering
        \includegraphics[width = \textwidth]{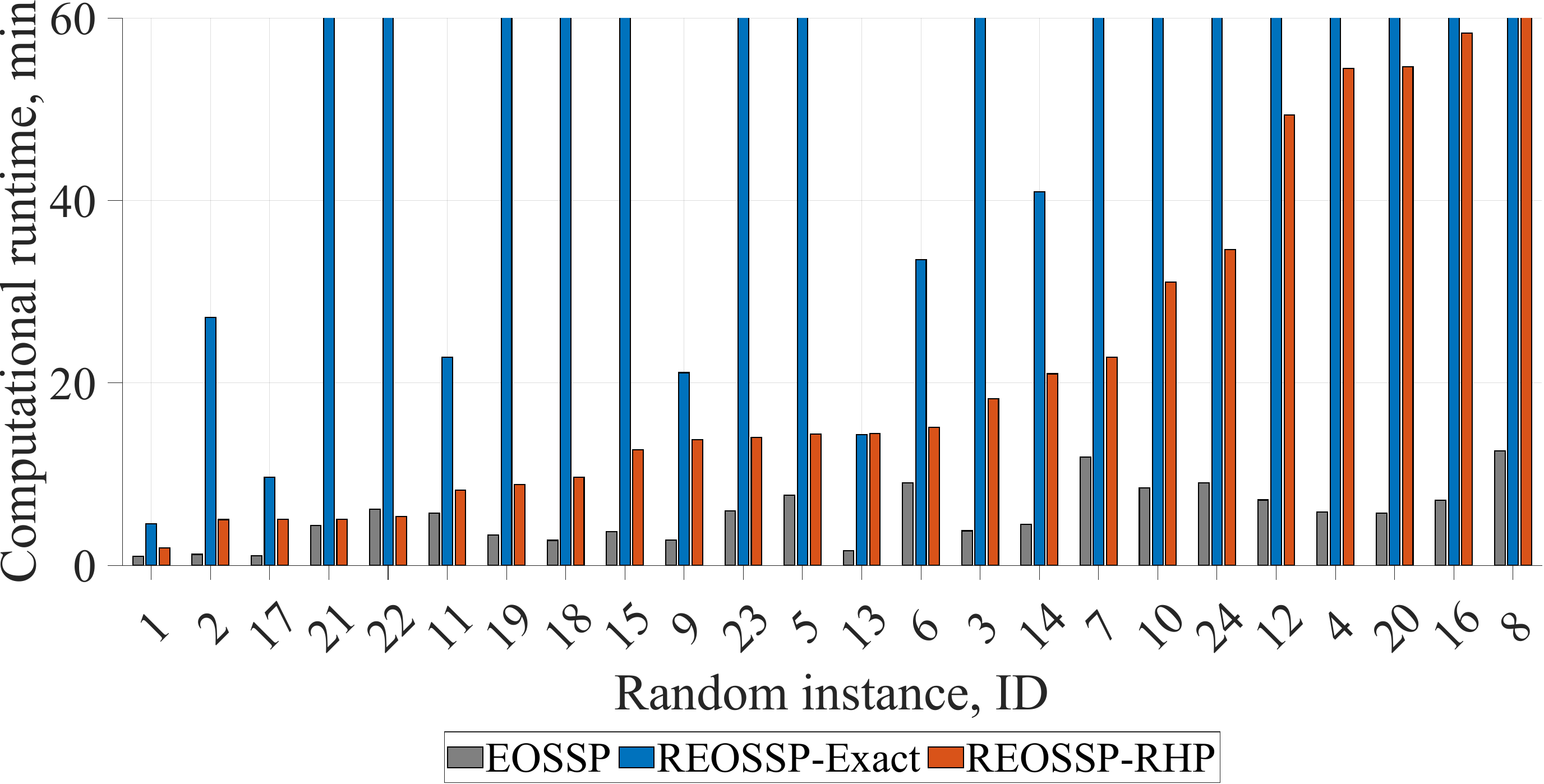}
        \caption{Runtime of the \EOSSP, \REOSSPExact, and \RHP}
        \label{fig:random_runtime}
    \end{subfigure}
    \caption{Results of the random instances.}
    \label{fig:rand_results}
\end{figure}

In addition to the statistical information provided by the box charts in Sec.~\ref{subsec:random_experiments}, the percent improvement of the \REOSSPExact and \RHP over the \EOSSP, as well as the runtime of the \EOSSP, \REOSSPExact, and \RHP, in each instance, is shown as bar charts in Figs.~\ref{fig:random_percent} and Fig.~\ref{fig:random_runtime}, respectively. The ordering of Figs.~\ref{fig:random_percent} and~\ref{fig:random_runtime} is in ascending order of the \RHP metric from left to right.

\section*{Acknowledgments}
The authors thank the anonymous reviewers for their contribution of feedback and insight, which has greatly improved the quality of the manuscript. The authors also thank the Department of Mechanical, Materials and Aerospace Engineering at West Virginia University for providing the support to publish this manuscript as open access.

\clearpage
\newpage
\bibliography{1_ProofReferences}

\end{document}